\documentclass[10pt]{amsart}
\usepackage[a4paper,margin=1in]{geometry}
\usepackage{amsmath,amssymb,amsthm,mathtools}
\usepackage{enumitem}
\usepackage[T1]{fontenc}
\usepackage{lmodern}
\usepackage{microtype}
\usepackage[hidelinks]{hyperref}
\usepackage{color}
\definecolor{oneblue}{rgb}{0,0,0.75}
\definecolor{Red}{cmyk}{0,1,1,0} 
\usepackage{young} \usepackage{ytableau}

\def\steq#1{ \stackrel{\eqref{#1}}}

\newcommand{\e}{\epsilon}
\newtheorem{theorem}{Theorem}[section]
\newtheorem{corollary}[theorem]{Corollary}
\newtheorem{proposition}[theorem]{Proposition}
\newtheorem{lemma}[theorem]{Lemma}

\theoremstyle{definition}
\newtheorem{definition}[theorem]{Definition}
\newtheorem*{example}{Example}
\theoremstyle{remark}
\newtheorem{remark}[theorem]{Remark}

\newcommand{\perm}{\textrm{perm}}

\newcommand{\Z}{\mathbb Z}
\newcommand{\N}{\mathbb N}
\newcommand{\Q}{\mathbb Q}

\newcommand{\GL}{\operatorname{GL}}
\newcommand{\SL}{\operatorname{SL}}
\newcommand{\Sp}{\operatorname{Sp}}
\newcommand{\OO}{\operatorname{O}}
\newcommand{\SO}{\operatorname{SO}}

\newcommand{\Div}{\operatorname{Div}}

\newcommand{\Span}{\operatorname{Span}}

\newcommand{\Ga}{\Gamma}

\newcommand{\ext}{\bigwedge}
\newcommand{\ot}{\otimes}
\title{Characteristic-Free Skew-Symmetric Classical Invariant Theory}
\author{Corrado De Concini 
Claudio Procesi}
\address{Dipartimento di Matematica ``Guido Castelnuovo'', Sapienza Universit\`a di Roma}

\date{}        
 
\begin{document}
\maketitle
 
\begin{abstract}
We study an exterior-algebra analogue of classical invariant theory for
linear, symplectic, and orthogonal groups, with particular emphasis on
integral constructions and positive characteristic.  

The theory presented here is rather different  from the classical one developed by us  in \cite{DeConciniProcesi1976}.

In fact, the present theory is based on good filtrations and on a surprising divided-power structure on   the algebras of exterior invariants.\medskip

 We also include a short characteristic-zero treatment of
the corresponding First and Second Fundamental Theorems via the exterior
Cauchy identity; these known results motivate the integral theory.
\end{abstract}

\medskip
\noindent\textbf{Keywords.} Exterior invariant theory; positive characteristic;
divided powers; integral standard monomials; good filtrations; Weyl
filtrations; classical groups.

\tableofcontents

\section*{Introduction}

The First and Second Fundamental Theorems of classical invariant theory
identify generators and relations for polynomial invariants of copies of the  defining representations of  classical
groups. Their integral formulations are well understood. In this paper we
study the corresponding questions when the symmetric algebra is replaced by
an exterior algebra, with particular emphasis on the integral structure and
on positive characteristic.\smallskip

We begin with the general linear group. For finite-dimensional vector spaces
$U,V,W$, with $\dim V=n$, consider
\[
 E(U,V,W)=\bigwedge\bigl(V\otimes U^*\oplus W\otimes V^*\bigr)
\]
and the natural $GL(V)$-action, commuting with that of
$GL(U)\times GL(W)$.\smallskip

 In matrix coordinates, Formula  \eqref{levar}, the basic quadratic contractions
are $t_{ij}=\sum_{s=1}^n y_{is}x_{sj}$. 

In characteristic zero the familiar
exterior analogues of the fundamental theorems follow by a short adaptation
of the classical methods; we recall this argument for comparison, without
claiming novelty.\smallskip

In positive characteristic the contractions alone need not generate the
invariants. Indeed, the commuting even elements $y_{is}x_{sj}$ have square
zero, and
\[
 t_{ij}^{\,r}=r!\,t_{ij}^{[r]},\qquad
 t_{ij}^{[r]}=\sum_{1\le s_1<\cdots<s_r\le n}
       \prod_{\nu=1}^{r}y_{is_\nu}x_{s_\nu j}.
\]
The appropriate characteristic-free formulation therefore uses divided
powers. A good-filtration argument, based on results of Wang \cite{Wang} (later greatly generalized by 
Donkin
and Mathieu; see~\cite{Donkin,Mathieu,Jantzen}), first shows that the dimensions of homogeneous
$GL(V)$-invariant spaces are independent of the characteristic.

Our main construction takes place over $\mathbb Z$. By dualizing the map from
the divided-power algebra on the matrix of contractions, we use the integral
standard bitableau basis of  polynomial  functions on matrices  to describe its
kernel and image. The admissible diagrams satisfy a rectangular bound. This
gives an integral basis of the image and a splitting in each degree; in
particular, the image is saturated. Combining this with the good-filtration
dimension calculation yields the invariant algebra, its divided-power
generators and its relations after base change to any field.

The special linear group case is considered as well. For the exterior algebra of vectors,
we identify the nonzero $\SL(V)$-invariant components and prove that they
are generated by the invariants in exterior degree $\dim V$. For mixed
vectors and covectors, the $\SL(V)$-invariant algebra is generated by
its $\GL(V)$-invariants together with  maximal skew determinants;
we also describe the resulting Weyl-module filtrations. 

Finally, we indicate the corresponding orthogonal and symplectic results.
These use the polynomial tableau theory for symmetric and alternating
matrices developed in earlier work of De Concini and Procesi. We also indicate how
the $\SO(V)$-invariants are obtained from $\OO(V)$-invariants by adjoining
 maximal skew determinants.
In these cases we state the
results and the method, omitting calculations parallel to the general linear
case which follow by analogous arguments.  In the case of the orthogonal group, characteristic two is excluded. Thus the
principal focus of the paper is the integral divided-power construction and
its characteristic-free consequences.  The authors acknowledge the use of  ChatGPT as a supporting tool for coherence checks, bibliographic searches, and expository revision.

  \medskip
\noindent
 
\part{Recollections}
\section{Classical invariant theory}

\subsection{The general linear group}

Let $K$ be a field and let $V\simeq K^n$ be an $n$-dimensional vector space.  We begin by recalling the standard characteristic-zero picture, which motivates the skew-symmetric version considered below.

Assume first that $\operatorname{char}K=0$.  For nonnegative integers $a=\dim U,b=\dim W$, 
consider \begin{equation}\label{lafond}
\mathcal A:=V^a\oplus V^{*b}=V\otimes U^*\oplus W\otimes V^*=\hom(U,V)\oplus \hom(V,W)\,.
\end{equation}   Choosing a basis for $V$, we  display the $a$ vectors in $V^a$  and the $b$   covectors in $V^{*b}$   as the columns $x_i$ of an $n\times a$ matrix $X$ and the rows  $y_j$ of a  $b\times n$ matrix $Y$:\begin{equation}\label{levar}
X=\{x_i=(x_{1,i}, x_{2,i},\dots, x_{n,i})\}; i=1,\cdots,a;\quad Y=\{y_j=(y_{j, 1},  y_{ j,2},\cdots, y_{ j,n})\};j=1,\cdots,b\,. 
\end{equation}  \medskip

The group $\mathbb G= \GL(V)\times \GL(U)\times \GL(W)$ acts naturally on $\mathcal A$ by the action on the two matrices $X,Y$:
\begin{equation}\label{matrix-action}
(A,B,C)(X)=A X B^{-1},\quad (A,B,C)(Y)=C Y A^{-1}
\end{equation} and hence on its coordinate ring. 

The map
\begin{equation}\label{acco}
\rho:(X,Y)\mapsto T:=Y\cdot X;\quad \rho:\hom(U,V)\times \hom(V,W)\to \hom(U,  W)
\end{equation}
is $\mathbb G$-equivariant where $ \GL(V)$ acts trivially on $\hom(U,  W)$. Let $T:=(t_{i,j})$ denote the resulting matrix.

\begin{theorem}[First Fundamental Theorem for $\GL(n)$]
The invariants $t_{i,j}$ generate the algebra of invariants $K[\mathcal A]^{ \GL(V)}$.\end{theorem}

Let $Z=(z_{i,j})$ be the generic $b\times a$ matrix.  The First Fundamental Theorem yields a surjection
\begin{equation}\label{acco1}
q:K[z_{i,j}]\longrightarrow K[\mathcal A]^{ \GL(V)},\qquad z_{i,j}\longmapsto t_{i,j}=\sum_{r=1}^ny_{i,r}x_{r,j}\,.
\end{equation}

\begin{theorem}[Second Fundamental Theorem for $\GL(n)$]
The kernel of $q$ is the ideal generated by the determinants of the $(n+1)\times(n+1)$ minors of $Z=(z_{i,j})$.
\end{theorem}
This theorem is classical, see for instance \cite{Weyl}  in characteristic 0, and was proved by the present authors in all characteristics in  \cite{DeConciniProcesi1976}.
\subsection{Orthogonal and symplectic groups}

Suppose now that $V$ carries a nondegenerate bilinear form $\langle-,-\rangle$, symmetric in the orthogonal case and skew-symmetric in the symplectic case.  In the latter case $n$ is even.   Denote by $G$ the corresponding group fixing the given bilinear form.\smallskip

Put $h=\dim U$. In these cases the bilinear form induces an isomorphism, under the orthogonal or symplectic group $G$,   between $V$ and $V^*$ and   we study the invariant theory for $ V^{\oplus h}
\simeq  \hom(U,V)$  and not $V\otimes U^*\oplus W\otimes V^*$.\smallskip

In this case  
we have the basic invariants
\[
 \rho:  X\mapsto T:=X^tJ X=(\langle v_i,v_j\rangle);\quad  \rho: \hom(U,V)\to \hom(U,U^*)\,.
\]
In the orthogonal case we take $J=1_n$ to be the identity matrix so $T=(t_{i,j})$ is a symmetric matrix.  In the  symplectic case the dimension $n=2k$ must be even and $J=\begin{vmatrix}
0&1_k\\-1_k&0
\end{vmatrix}$  is the standard     skew-symmetric matrix, and $T$ is also skew-symmetric.  The map $\rho$ is also $\GL(U)$-equivariant.

\begin{theorem}
The functions $t_{i,j}$ generate the invariant ring $A_h=K[V^{\oplus h}]^G$, where $G=\OO(V)$ or $\Sp(V)$ respectively.
\end{theorem}

\begin{theorem}[SFT for $\OO(n)$ and $\Sp(n)$]\label{FOG}
Let $Z=(z_{i,j})$ denote the corresponding generic symmetric or skew-symmetric matrix.
\begin{enumerate}[label=\textnormal{(\roman*)}]
\item For $\OO(n)$, the ideal of relations is generated by the determinants of the minors of size $n+1$ of $Z$.
\item For $\Sp(n)$, the ideal of relations is generated by the Pfaffians of the principal minors of size $n+2$ of $Z$.
\end{enumerate}
\end{theorem}

These results belong to the classical invariant theory of the classical groups; see, in particular, Weyl's \emph{The Classical Groups} \cite{Weyl} (Part  (ii) of Theorem \ref{FOG} is not proved there). \smallskip

 A characteristic-free treatment of the First and Second Fundamental Theorems over a  field, and in integral form, was developed in \cite{DeConciniProcesi1976}.  The orthogonal case in characteristic $2$ requires special care, and indeed is not covered in \cite{DeConciniProcesi1976} (contrary to what is stated there).

\subsection{Aronhold's method}  A method widely used in classical invariant theory, in characteristic 0,   is  polarization and restitution,  which reduces the study of invariants  and their relations  to the corresponding multilinear problem.
\smallskip

Given a linear representation $V$  of an algebraic group $G$, one wants to study the invariants of the polynomial algebra $S[V^*]$. The idea is to study $G$-- invariants of homogeneous degree $d$ by polarization to multilinear invariants of $d$  copies of $V$. Once one has a description of such invariants the homogeneous invariants are obtained by the process of {\em restitution}.\smallskip

The method holds verbatim  if one replaces $S[V^*]$ with $\bigwedge  V  $ under some action of a group $G$ on $V$.  

Here the polarization map for elements of degree $d$  is obtained by considering the map 
$$P: \bigwedge   V  \to \bigwedge   V^{\oplus d}  [\lambda_1,\cdots,\lambda_d];\  v\mapsto \sum_{i=1}^d\lambda_i v_i$$ where the $\lambda_i$ are parameters.

The {\em restitution} is the map  
$$R: \bigwedge   V^{\oplus d}  \to \bigwedge   V  ,\quad (v_1,v_2,\cdots, v_d)\mapsto\sum_iv_i\,.$$
Given an element  $a\in \bigwedge^dV$ invariant  under $G$, the element $P(a)$ is a polynomial in the $\lambda_i$ with coefficients in $\bigwedge   V^{\oplus d}  $ which are also $G$-invariant, in particular the coefficient of $\prod_{i=1}^d\lambda_i$  is the full polarization $\bar P$ and it is a multilinear invariant. One has $R\bar P(a)=d!a$.\medskip

The method is then the following: one describes a space of invariants closed under polarization and restitution  and then if the multilinear elements of this space coincide with the space of all multilinear invariants, the proposed invariants coincide with the full space of invariants.
\smallskip

One has the same space $(V^{\otimes d})^G$ of multilinear invariants for the exterior  or the symmetric case and the only problem is to {\em decode} the formulas  passing from symmetric to exterior formulas.  \smallskip

This can be done by a simple trick: replace the $d$ vector Grassmann variables $v_i$  with {\em symmetric variables $\epsilon_iv_i$}  where the $\epsilon_i$  are auxiliary Grassmann variables.  
\smallskip

We quickly illustrate this idea  for the linear group (see Theorem \ref{thm:skew-fft-sft} for details). 

 In this case first notice that the linear group contains the group $K^*$ of scalars $\alpha$ which acts on $V$ by multiplication and on $V^*$   by multiplication with  $\alpha^{-1}$ so a homogeneous invariant must have the same degree in $V$ and $V^*$. Polarization must then be performed independently on $V$ and $V^*$.
 
   One has the multilinear invariants  in $d$ copies of $V$ and $d$ copies of $V^*$  which can be identified with $End_{GL(V)}V^{\otimes  d}$  (or its dual) which is spanned by elements of the symmetric group  $S_d$   acting on $V^{\otimes  d}$.

 The kernel of the map from the algebra of the symmetric group to the linear operators on $V^{\otimes  d}$  is the ideal generated by the antisymmetrizer.\smallskip

  In the symmetric case the invariant associated to a permutation $\sigma$  is $t_{1,\sigma(1)}t_{2,\sigma(2)}\cdots t_{d,\sigma(d)} $  as follows from  formula \eqref{acco} and   in the Grassmann case this acquires the sign $\epsilon_\sigma$ of the permutation.
  \smallskip

  This produces the first and second fundamental theorems in the form which will be discussed later in this article. The two different interpretations of the antisymmetrizer give rise to the two different formulations of the SFT.
  \smallskip

  The presentation in Section \ref{c0i} is more precise since it also describes the action of the group $GL(U)\times GL(W)$ on the invariants of $GL(V)$.
  
  \section{Exterior Cauchy identities and skew invariant theory\label{c0i}}

We now replace the symmetric algebra of polynomial functions by an exterior algebra.  Throughout this section $K$ has characteristic zero.  

For a partition $\lambda=(\lambda_1\geq\lambda_2\geq\cdots)$, write $\lambda'$ for its transpose, let $\ell(\lambda)=\lambda'_1$ be the height of the partition.  We use the convention that $s_\lambda$ denotes the Schur function indexed by the transpose of $\lambda$ in the usual convention (so $s_{(r)}=e_r$), and that $L_\lambda(E)$ is the irreducible polynomial $\GL(E)$-module whose character is $s_\lambda$.  Thus
\[
L_{(r)}(E)=\bigwedge^r E,
\qquad
L_{(1^r)}(E)=S^rE.
\]
With this convention, the exterior Cauchy identity reads (see, for example, \cite[Chapter~4]{Fulton})
\begin{equation}\label{Cauchy}
\prod_{i=1}^{n}\prod_{j=1}^{h}(1+x_i y_j)
=
\sum_{\substack{\lambda:\ \ell(\lambda)\le h,\ \lambda_1= \ell(\lambda')\le n}}
 s_\lambda(x_1,\ldots,x_n)s_{\lambda'}(y_1,\ldots,y_h),
\end{equation}
and gives the decomposition into irreducible representations of $GL(V)\times GL(U)$.
\begin{equation}\label{eq:exterior-cauchy}
\bigwedge^*(V\otimes U^*)
\simeq
\bigoplus_{\substack{\lambda:\  \ell(\lambda')\le \dim V,\ \ell(\lambda)\le \dim U}}
L_\lambda(V)\otimes L_{\lambda'}(U^*).
\end{equation}Here $\lambda_1\leq\dim V=n$ and $\ell(\lambda)\leq\dim U=a$: the diagram of $\lambda$ lies in the rectangle of width $n$ and height $a$.\medskip

Recall that  $L_\lambda(V)\otimes L_{\mu}(V^*)=\hom(L_\mu(V),L_\lambda(V))$ so $(L_\lambda(V)\otimes L_{\mu}(V^*))^{\GL(V)}$ is 0 if $\mu\neq \lambda$  while $L_\lambda(V)\otimes L_{\lambda}(V^*)=End(L_\lambda(V) )$  so its invariants are the 1--dimensional space generated by the identity.\medskip

We now use the notations of formula \eqref{lafond} and set:
\begin{equation}\label{IE1}
  \bigwedge \mathcal A= \bigwedge(V\otimes U^*\oplus W\otimes V^*)= \bigwedge(\hom(U,V)\oplus \hom(V,W))=: E(U,V,W)\,.
 \end{equation}

Applying \eqref{eq:exterior-cauchy} to the two summands of the space $E(U,V,W)$ of formula \eqref{IE1} and taking $\GL(V)$-invariants gives, 
$$\bigwedge(V\otimes U^*\oplus W\otimes V^*)=\bigwedge(V\otimes U^*)\otimes \bigwedge(W\otimes V^*)\simeq\bigoplus  L_{\lambda'}(U^*)\otimes L_\lambda(V)\otimes  L_\mu(V^*)\otimes L_{\mu'}(W)  $$
\begin{align}\implies\qquad
E(U,V,W)^{\GL(V)}
&\simeq
\bigoplus_{\substack{\lambda_1\le \min(a,b)\\ \ell(\lambda)\le n}}
L_\lambda(U^*)\otimes L_\lambda(W ).
\label{eq:invariant-decomposition}
\end{align} Here $\lambda_1\leq\min(a,b)$ and $\ell(\lambda)\leq n$: the diagram of $\lambda$ lies in the rectangle of width $\min(a,b)$ and height $n$.

\subsection{The basic quadratic invariants\label{bqi}}With the notation of formula \eqref{lafond} we see that the bilinear map of formula   \eqref{acco} can be interpreted also  as a map in the exterior algebra, defining now the invariant elements  

 $t_{i,j}  =\sum_{r=1}^ny_{i,r}x_{r,j}$ of  even degree in the exterior algebra. \smallskip
 
 Thus  $t_{i,j} $ is central and in particular the $t_{i,j}$ mutually commute.  They span the copy of $U^*\otimes W$ occurring in degree two in \eqref{eq:invariant-decomposition}.  
 
 Hence there is a natural $\GL(U)\times\GL(W)$-equivariant homomorphism as in formula \eqref{acco1}:
\begin{equation}\label{ilq}
q:K[z_{i,j}]\longrightarrow E(U,V,W)^{\GL(V)},
\qquad z_{i,j}\longmapsto t_{i,j}=\sum_{r=1}^ny_{i,r}x_{r,j}\,.
\end{equation}

\begin{theorem}
\label{thm:skew-fft-sft}
With the notation above:
\begin{enumerate}[label=\textnormal{(\roman*)}]
\item The homomorphism $q$ is surjective.  Equivalently, the elements $t_{i,j}$ generate $E(U,V,W)^{\GL(V)}$.
\item The kernel of $q$ is the ideal generated by the irreducible $\GL(U)\times\GL(W)$-submodule
\[
S^{n+1}(U^*)\otimes S^{n+1}(W)
=
L_{(1^{n+1})}(U^*)\otimes L_{(1^{n+1})}(W).
\]
\end{enumerate}
\end{theorem}

\begin{proof} i)\quad The proof follows from Aronhold's method. We use a general trick which applies also to the orthogonal and the symplectic groups.\smallskip

We want to compute the multilinear invariants in $\mathcal A =\bigwedge (V^{\oplus d}\bigoplus V^{* \oplus d})$.  \smallskip

Add Grassmann variables  $\epsilon_1,\cdots,\epsilon_d;\eta_1,\cdots,\eta_d$  and then in $\mathcal A \otimes \bigwedge[\epsilon_1,\cdots,\epsilon_d;\eta_1,\cdots,\eta_d]$ consider the vectors $\bar v_i:=\epsilon_iv_i,\ \bar \phi_i:=\eta_i \phi_i$, these are all even, so we have a homomorphism $\Sigma=S[V^{\oplus d}\bigoplus V^{* \oplus d}]\to \mathcal A \otimes \bigwedge[\epsilon_1,\cdots,\epsilon_d;\eta_1,\cdots,\eta_d]$  mapping $v_i\to \bar v_i,\ \phi_i\to\bar\phi_i$. Under this map  the multilinear space $V^{\otimes d}\otimes V^{* \otimes d}\subset\Sigma$ maps isomorphically to  the multilinear space $V^{\otimes d}\otimes V^{* \otimes d}\prod_{i=1}^d\epsilon_i\eta_i\subset \mathcal A \otimes \bigwedge[\epsilon_1,\cdots,\epsilon_d;\eta_1,\cdots,\eta_d]$. \smallskip

Under this map an element  $t_{i,j}\in \Sigma$   maps to $-t_{i,j}\eta_i\epsilon_j\in  \mathcal A \otimes \bigwedge[\epsilon_1,\cdots,\epsilon_d;\eta_1,\cdots,\eta_d]$.\smallskip

  This map is $\GL(V)$-equivariant so it is an isomorphism on the invariants and maps a multilinear product  $u_\sigma=\prod_it_{i,\sigma(i)},\ \sigma\in S_d$ of the elements $t_{i,j}\in \Sigma$  to the   product in $\mathcal A $  but multiplied by $\prod_{i=1}^d(-\eta_i\epsilon_{\sigma(i)})=(-1)^{d(d+1)/2}\epsilon_\sigma \eta_1\cdots\eta_d\epsilon_1\cdots\epsilon_d$, with $\epsilon_\sigma$  the sign of the permutation.\smallskip

This shows that the multilinear invariants  in $\mathcal A $  are spanned by the $u_\sigma=\prod_it_{i,\sigma(i)},\ t_{i,j}\in \mathcal A $ then by Aronhold's method this implies that the  $t_{i,j}\in \mathcal A $ generate all invariants.\medskip

ii)\quad For the relations, one can also use Aronhold's method  by noticing that the multilinear relations in the various degrees are of the form  $\sum_\sigma c_\sigma u_\sigma$ where in the symmetric case the element $\sum_\sigma c_\sigma  \sigma$ is in the ideal of the group algebra generated by the antisymmetrizer on $n+1$  indices, $n=\dim V$; while in the exterior case  $\sum_\sigma c_\sigma  \sigma$ is in the ideal of the group algebra generated by the  symmetrizer on $n+1$  indices.\medskip

On the other hand  one can argue by noticing that
\[
q(z_{1,1}^{n+1})=t_{1,1}^{n+1}
 =\left(\sum_{i=1}^n y_{1,i}x_{i,1}\right)^{n+1}=0.
\]
The element $z_{1,1}^{n+1}$ is a highest-weight vector of
\[
S^{n+1}(U^*)\otimes S^{n+1}(W)
=L_{(1^{n+1})}(U^*)\otimes L_{(1^{n+1})}(W ).
\]
By Theorem 4.1 of \cite{DePE1980}, based on Pieri's rule, the ideal generated by this module in $K[z_{i,j}]$ contains precisely
the irreducible summands indexed by partitions with at least $n+1$
rows, i.e. the summands excluded from
\eqref{eq:invariant-decomposition}; see also
\cite[Chapter~5]{Fulton}.  \smallskip

Comparing the ordinary Cauchy decomposition
of $K[U^*\otimes W ]$ with \eqref{eq:invariant-decomposition} proves the
claim.
\end{proof}

\section{The orthogonal and symplectic cases}
\label{sec:char-zero-orth-symp}

The same method applies to the other classical groups.  We use the
following standard branching facts, in the convention for $L_\lambda$
adopted above; see Weyl~\cite[Chapters~VI--VII]{Weyl}.

If $V$ is equipped with a nondegenerate symmetric bilinear form, then
\begin{equation}\label{InvOrt}
L_\lambda(V)^{\OO(V)}\ne0
\quad\Longleftrightarrow\quad
\text{all columns of $\lambda$ have even length},
\end{equation}
and in this case the invariant space is one-dimensional.  If $V$ is
equipped with a nondegenerate alternating form, then
\begin{equation}\label{InvSym}
L_\lambda(V)^{\Sp(V)}\ne0
\quad\Longleftrightarrow\quad
\text{all rows of $\lambda$ have even length},
\end{equation}
and again the invariant space is one-dimensional.

\begin{theorem}
\label{thm:orthogonal}
Let $K$ have characteristic zero and let $V$ carry a nondegenerate
symmetric bilinear form. Then
\[
\left(\bigwedge^*(V\otimes W)\right)^{\OO(V)}\simeq
\bigoplus_{\substack{\mu\ \mathrm{row\ even}\\
\ell(\mu)\leq\dim V,\ \mu_1\leq\dim W}}L_\mu(W)
\simeq K[\bigwedge^2W]/I,
\]
where $I$ is generated by the simple $\GL(W)$-module
\[
L_{(2^{\,\dim V+1})}(W).
\]
\end{theorem}

\begin{proof}
The exterior Cauchy decomposition \eqref{eq:exterior-cauchy}
 and \eqref{InvOrt} give
\[
\left(\bigwedge^*(V\otimes W)\right)^{\OO(V)}
\simeq
\bigoplus_{\substack{\mu\ \mathrm{row\ even}\\
\ell(\mu)\leq\dim V,\ \mu_1\leq\dim W}}L_\mu(W).
\]
Choose a basis of $V$ with Gram matrix $J$ and a basis of $W$, and write
$X=(x_{i,j})$ for the matrix of exterior generators as in page \pageref{lafond}.  Then
\[
T={}^tXJX=(t_{r,s})
\]
is skew-symmetric, and its entries $t_{r,s}$, $r<s$, span
$\bigwedge^2W$.  Hence they define a homomorphism
\[
q:K[\bigwedge^2W]\longrightarrow
\left(\bigwedge^*(V\otimes W)\right)^{\OO(V)}.
\]
The same trick of adding Grassmann variables $\epsilon_i$  and replacing $v_i$ with $\epsilon_iv_i$  proves
surjectivity exactly as in Theorem~\ref{thm:skew-fft-sft}.\medskip

Finally $t_{1,2}^{\dim V+1}=0$.  Its preimage is a highest-weight vector
of $L_{(2^{\dim V+1})}(W)$, and Pieri's rule, cf.  \cite{Abd} Theorem 3.1, shows that the ideal it
generates contains exactly the excluded summands.  \smallskip

This proves the
statement.  Compare also the classical characteristic-free treatment
in \cite{DeConciniProcesi1976}.
\end{proof}

\begin{theorem}
\label{thm:symplectic}
Let $K$ have characteristic zero and let $V$ carry a nondegenerate
alternating form. Then
\[
\left(\bigwedge^*(V\otimes W)\right)^{\Sp(V)}
\simeq K[S^2W]/I,
\]
where $I$ is generated by
\[
L_{(1^{\,\dim V+2})}(W)=S^{\dim V+2}(W).
\]
\end{theorem}

\begin{proof}
The proof is parallel to the orthogonal case, using \eqref{InvSym}.
If $J$ is a matrix of the symplectic form, then
\[
T={}^tXJX
\]
is symmetric because the entries of $X$ anticommute; its entries span
$S^2W$.  When applying Aronhold's method, note that  a quadratic invariant $x_i^tJx_i $  is not  bilinear but its polarization is nonzero. \smallskip

 The first excluded representation is
$L_{(1^{\dim V+2})}(W)$, and Pieri's rule, cf. \cite{Ab}  gives the ideal of relations.
See also \cite{Weyl,DeConciniProcesi1976}.
\end{proof}

\part{Positive characteristic and good filtrations}\label{pogf}
When we deal with invariant theory for an affine  variety $V$, defined over a  field $K$ with coordinate ring $K[V]$,    under an action of an algebraic group $G$ defined over  $K$ we do not usually want to consider the invariants under the group $G(K)$ of its rational points but under $G(\bar K)$ with $\bar K$ the algebraic closure. These invariants are the same if  $G(K)$ is Zariski dense in  $G(\bar K)$. \smallskip

For classical groups  this is true  as soon as $K$ is an infinite field, but in our examples the groups and the varieties, which are just linear spaces, are defined over $\Z$ and hence over all fields and 
so also  over the prime  field $\Z/(p)$. 

  The (group scheme) action is defined by the map $\delta: \Z[V]\to \Z[V]\otimes \Z[G]$  which induces for every field $K$ an action $\delta_K: K[V]\to K[V]\otimes K[G]$  and an element   $a\in \Z[V]$ is an invariant if and only if  $\delta(a)=a\otimes 1$.

This in our case is equivalent to the fact that for every infinite field  $H\supset K$  the element $a$ is invariant  under the rational points $G(H)$.\smallskip

In the orthogonal case, these arguments require characteristic different from two.

 In fact the coordinate ring of the  group $\OO(n)$ over $\Z$  is the quotient of the polynomial ring $\Z[x_{i,j}],\ i,j=1,\cdots,n$ modulo the quadratic equations  $X\cdot X^t=1$.   
 
 Passing to   characteristic 2, this ring is not even reduced  since one has for all $i=1,\cdots,n$  that $(\sum_jx_{i,j}+1)^2=0$.\medskip
 
 Apart from this pathological case all other groups are reductive in all characteristics.

See \cite{dem},\  \cite{Co},  and  \cite{Co1}, \S 5 on classical groups, which treats $Sp_{2n}$, $SL_n$ and $SO(q)$ (with 2 invertible for $SO(q)$).

\label{SchurWeyl}\section{Schur and Weyl modules}
For the reader's convenience,  we  recall some standard facts concerning Weyl and Schur modules for $GL(V)$   emphasizing their combinatorial aspects. A standard reference for this is  \cite{ABW}.

\subsubsection{The divided-power algebra}  

Given a commutative $\Q$-algebra $A$ and an element $a\in A,\ r\in \N$ one defines the {\em divided power $a^{[r]}:=\frac {a^r}{r!}$}. 

 One has the binomial expansion $(a+b)^{[r]}:=\sum_{i=0}^r a^{[r-i]}b^{[i]}$.\medskip

The divided-power algebra in $n$ variables is a Hopf algebra and  can be defined in several ways. We define it over $\Z$ and denote it by  $\Ga \Z^n=\oplus_{d=0}^\infty \Ga^d \Z^n$. It is naturally acted upon as automorphisms by the group scheme $GL(n,\Z)$.
\begin{definition}\label{dpa}
\begin{enumerate}\item  $\Ga \Z^n$ is the $\Z$  subalgebra of the polynomial algebra $\Q[e_1,\cdots,e_n]$  generated by all the divided powers of elements of $\bigoplus_{i=1}^n\Z e_i$.\item   $\Ga \Z^n$ is naturally identified with the graded dual of  $\Z[e^1,\cdots,e^n]$  where $e^i$ is the dual basis of $e_i$.

\item The space $\Gamma^d(\Z^n)$ is isomorphic to the symmetric tensors in  $(\Z^n)^{\otimes d}$. Under this isomorphism,  multiplication is by symmetrization of the tensor product.\end{enumerate}\end{definition}

In fact  these different definitions  stem from the fact that $\Ga \Z^n$ is a Hopf algebra with the structure induced by the Hopf structure of  $\Q[e_1,\cdots,e_n]$ for Definition 1) and by duality   for Definition 2). 

The comultiplication $\Delta$  is given by   $$\Delta(a ^{[r]}):=\sum_{i=0}^r a^{[r-i]}\ot a^{[i]}.$$

For   3) the isomorphism is given by  the classical {\em full polarization}, which in the language of Hopf algebras consists in taking  the multilinear part of multidegree $(1,\ldots,1)$  denoted by $\Delta^d_m$ of the iterated  comultiplication $\Delta^d$.

 The product becomes   the symmetrized tensor product, where if $A\in \Ga^a \Z^n,\ B\in \Ga^b \Z^n$  the symmetrized product
$$A\cdot B\in \Ga^{a+b}\Z^n;\quad A\cdot B=\sum_{\sigma\in S_{a+b}/S_{a }\times S_{ b}} \sigma(A\ot B)$$ the sum over coset representatives.
\subsubsection{A combinatorial  approach}  In order to develop a characteristic-independent formalism   
it is important to work over $\Z$.  \medskip

So,  for  any partition $\lambda=(\lambda_1\geq\lambda_2\geq\cdots)$ we are going to   define two integral forms $ \nabla_\lambda(\Z^n),\ \Delta_\lambda(\Z^n)$,  of the irreducible $L_\lambda(V)$
which  are free abelian groups equipped with an action of the group scheme $GL(n,\Z)$. 

 Then for every field $K$  one has $\nabla_\lambda(K^n)=\nabla_\lambda(\Z^n)\otimes_\Z K;\quad \Delta_\lambda(K^n)=\Delta_\lambda(\Z^n)\otimes_\Z K$.\medskip

If $\lambda=(d)$ is a partition with a single row  we have  $\nabla_d(\Z^n)=\Delta_d(\Z^n)=\bigwedge^d(\Z^n)$ while $\nabla_{1^d}(\Z^n)=S^d(\Z^n)$ the symmetric power.  Finally $\Delta_{1^d}(\Z^n)=\Gamma^d(\Z^n)$ the symmetric tensors in  $(\Z^n)^{\otimes d}$. \medskip

For general $\lambda$  one has to use an integral form of the classical {\em Young symmetrizers}.
One starts with a Young tableau ${\mathtt T}$, that is a Young diagram of some shape $\lambda\vdash d$, this decomposes the set $\{1,2,\cdots,d\}$  into two partitions  the {\em rows}  and the {\em columns} of the tableau.   Then  one defines the two subgroups $C_{\mathtt T},\ R_{\mathtt T}$ of the symmetric group $S_d$  preserving respectively the columns and the rows  and the elements in the group algebra:

\begin{equation}\label{youngs}
a_{\mathtt T}=\sum_{ \sigma\in R_{\mathtt T}}\e_\sigma \sigma,\ \text{antisymmetrizer};\  s_{\mathtt T}=\sum_{ \sigma\in C_{\mathtt T}} \sigma,\ \text{symmetrizer.}
\end{equation} 

In the classical theory one builds from ${\mathtt T}$ a {\em Young symmetrizer}  $c_{\mathtt T}$  and  then $c_{\mathtt T}V^{\otimes d}  =L_\lambda(V)$  is the corresponding $GL(V)$ irreducible, embedded in $V^{\otimes d}$ in a way which depends on ${\mathtt T}$.

There are in fact   two different primitive idempotents which can be called   $c_{\mathtt T}$.   One is up to scalar  $a_{\mathtt T}s_{\mathtt T}$  and the other $s_{\mathtt T}a_{\mathtt T}$.
\medskip

  Using these two partitions one defines (with a slight abuse of notation) two maps
\begin{equation}\label{duef}
\pi_{\mathtt T}:\bigotimes_j\bigwedge^{\lambda_j}\Z^n\stackrel{a_{\mathtt T}}\to (\Z^n)^{\otimes d}\stackrel{s_{\mathtt T}}\to \bigotimes_jS^{\lambda'_j}\Z^n;\quad \delta_{\mathtt T}:\bigotimes_j\Gamma^{\lambda'_j}\Z^n\stackrel{s_{\mathtt T}}\to (\Z^n)^{\otimes d}\stackrel{  a_{\mathtt T}}\to \bigotimes_j\bigwedge^{\lambda_j}\Z^n\,
\end{equation}
The arrows in \eqref{duef} are the canonical embeddings and quotient maps. In the first $\pi_{\mathtt T}$  the map $a_{\mathtt T}$  identifies $\bigotimes_j\bigwedge^{\lambda_j}\Z^n$ with the tensors in $(\Z^n)^{\otimes d}$  which are antisymmetric in the indices of each row  while $s_{\mathtt T}$  projects  to the symmetric quotient in the indices of each column.

Conversely  for $\delta_{\mathtt T}$  the map $s_{\mathtt T}$  identifies $\bigotimes_j\Gamma^{\lambda'_j}\Z^n$ with the tensors in $(\Z^n)^{\otimes d}$  which are  symmetric in the indices of each column  while $a_{\mathtt T}$  projects  to the exterior quotient in the indices of each row.
\smallskip

Unlike exterior powers which can be identified with submodules or quotient modules of  tensor powers,  $\Gamma^{a}\Z^n$ and $S^a\Z^n$  are different over the integers: the former arise as submodules  and the latter  as quotients.   
\begin{definition}\label{nabdel}
We define  $\nabla_\lambda(\Z^n)$ to be the image of  $\pi_T$ and  $\Delta_\lambda(\Z^n)$ to be the image of  $\delta_T$.
\end{definition}
Notice that  these modules depend only on the shape $\lambda$ and not on the choice of the tableau $T$.

This is the basis for the  combinatorial description   of the modules  $\nabla_\lambda(\Z^n),\ \Delta_\lambda(\Z^n)$ via standard tableaux which we briefly describe.

 \begin{remark}\label{duas}
Both maps in \eqref{duef} can be dualized  giving rise to an isomorphism of modules $\Delta_{\lambda}(\Z^n)\simeq \nabla_{\lambda}((\Z^n)^*)^*.$
\end{remark}

\medskip

 i)\quad {\bf the module $\nabla_\lambda(\Z^n)$} \quad A basis  of  $\bigotimes_j\bigwedge^{\lambda_j}\Z^n$  is formed by the elements   $$\bigotimes_je_{i_1}\wedge e_{i_2} \wedge\cdots \wedge e_{i_{\lambda'_j}},\ 1\leq i_1 <{i_2} <\cdots <{i_{\lambda'_j}}\leq n$$  which one displays as a  tableau $T$  of  shape $\lambda$ as in the following example.
$$\raisebox{ 15pt}{$T=$}\quad\begin{Young}
2&3&4&5 \cr
1&2&4&7\cr
4&6\cr \end{Young} \raisebox{ 15pt}{ ;}\qquad \raisebox{ 15pt}{$e_2\wedge e_3\wedge e_4\wedge e_5\otimes e_1\wedge e_2\wedge e_4\wedge e_7\otimes e_4\wedge e_6$}$$
Among these tableaux there are those called {\em semistandard} or sometimes just  standard. These are the tableaux in which the rows are strictly increasing while the columns are only weakly increasing. For instance, the preceding tableau $T$ is not semistandard while the next is:
$$\raisebox{ 15pt}{$\bar T=$}\quad\begin{Young}
2&3&4&5 \cr
2&3&4&7\cr 4&6\cr \end{Young} \raisebox{ 15pt}{ ;}\qquad \raisebox{ 15pt}{$e_2\wedge e_3\wedge e_4\wedge e_5\otimes e_2\wedge e_3\wedge e_4\wedge e_7\otimes e_4\wedge e_6$}$$
ii)\quad {\bf the module $\Delta_\lambda(\Z^n)$} \quad   A basis  of  $\Gamma^a\Z^n$  can be described by writing
\begin{equation}\label{lepp}
\Delta^d_m(\sum_{j=1}^n\alpha_je_j)^{[ a]}=(\sum_{j=1}^n\alpha_je_j)^{\otimes a}=\sum_{h_1,\cdots,h_n\mid\sum h_i=a}\prod_j\alpha_j^{h_j}[\underbrace{1,\cdots,1}_{h_1},\underbrace{2,\cdots,2}_{h_2},\cdots, \underbrace{n,\cdots,n}_{h_n}]
\end{equation}
where the underbraced symbols represent the corresponding symmetric tensors, an integral basis of  symmetric tensors.

In the Hopf algebra language  the    underbraced symbols  are $\Delta^d_m(e_1^{[h_1]} e_2^{[h_2]} \cdots e_n^{[h_n]} )$ or the symmetrized product $$[\underbrace{1,\cdots,1}_{h_1},\underbrace{2,\cdots,2}_{h_2},\cdots, \underbrace{n,\cdots,n}_{h_n}]=[\underbrace{1,\cdots,1}_{h_1}]\cdot[\underbrace{2,\cdots,2}_{h_2}]\cdot \cdots \cdot[ \underbrace{n,\cdots,n}_{h_n} ]\,.$$\medskip

Again we may display an element of the induced basis of  $\bigotimes_j\Gamma^{\lambda'_j}\Z^n$  as a tableau of shape  $\lambda'$  with rows  the symbols $[ i_1 ,{i_2} ,\cdots ,{i_{\lambda'_j}} ];1\leq i_1 \leq {i_2} \leq \cdots \leq {i_{\lambda'_j}}\leq n$.
 In this case the rows are weakly increasing.
The fundamental fact of the theory of standard tableaux is
\begin{theorem}\label{stat}
1) The images in  $\nabla_\lambda(\Z^n)$ of the semistandard  tableaux  form an integral basis of $\nabla_\lambda(\Z^n)$.

2) The images in  $\Delta_\lambda(\Z^n)$ of the   tableaux $T$  such  that their transpose is semistandard form an integral basis of $\Delta_\lambda(\Z^n)$.
\end{theorem}

Spanning follows by induction from the following proposition for the maps $\pi_T$ and $\delta_T$, which describes the {\em straightening relations}  which produce an explicit algorithm to express the image of a nonsemistandard  tableau as a linear combination of semistandard  tableaux. One orders the tableaux by lexicographic order of their {\em row words}  for instance, the  row word of $T$  is $2,3,4,5,1,2,4,7,4,6$.
\smallskip

 The expansion is unique because over $\Q$ the number of standard tableaux of a given shape equals the dimension of $\nabla_\lambda$ (resp. $\Delta_\lambda$). Linear independence over $\Q$ implies integral linear independence.
\begin{proposition}\label{eze} 1) Let $I=(i_1<\cdots< i_a)$,
$J=(j_1<\cdots< j_b)$, with $a\ge b$. Assume that the tableau $[I][J]=[i_1\ldots i_a][j_1\cdots j_b]$ is not standard. \smallskip

Then,  in $\nabla_{a,b}(\Z^n)$,
\begin{equation}\label{eq:tri-schur}
 [I][J]=\sum_{[I',J']\ne[I,J]}c_{I',J'}[I'][J']
 \qquad c_{I',J'}\in\mathbb{Z},
\end{equation}  where every $I'$ appearing in the sum satisfies $I'<_{\mathrm{lex}}I$.\medskip

2) Let $I=(i_1\le\cdots\le i_a)$,
$J=(j_1\le\cdots\le j_b)$, with $a\ge b$. Assume that the transpose of the tableau $[I][J]=[i_1\ldots i_a][j_1\cdots j_b]$ is not standard.  \smallskip

Then,  in $\Delta_{2^b1^{a-b}}(\Z^n)$,
\begin{equation}\label{eq:tri}
 [I][J]=\sum_{[I',J']\ne[I,J]}c_{I',J'}[I'][J']
 \qquad c_{I',J'}\in\mathbb{Z},
\end{equation}  where every $I'$ appearing in the sum satisfies $I'<_{\mathrm{lex}}I$.
\end{proposition}
\begin{proof} We give the proof of the second assertion since the proof of the first is  analogous and simpler.
Let $h$ be the smallest index such that $j_h\leq i_h$. Set
\[
 t=\max\{s\in\{h,\ldots,b\}:j_s=j_h\}.
\]
Let $A,B,C,D$ be the divided-power monomials associated with the multisets
\[
 (i_1,\ldots,i_{h-1}),\quad(i_h,\ldots,i_a),\quad
 (j_1,\ldots,j_t),\quad(j_{t+1},\ldots,j_b),
\]
respectively, and let $E$ be the monomial associated with the union of the multisets corresponding to $B$ and $C$.
Take the component 
\begin{equation}\label{lasomma}\mathcal S=\sum_\ell E_\ell^{(1)}\ot E_\ell^{(2)}\end{equation}
of the coproduct $\Delta(E)$  of bidegree $(a-h+1,t)$. 
Notice that  exactly one summand of \eqref{lasomma}, say $E_0^{(1)}\ot E_0^{(2)}$ equals $B \otimes C$ while for $\ell\geq 1$, in $E_\ell^{(1)}$ at least one index in $B$ has been replaced by a strictly smaller one. The choice of $h$ and $t$ ensures that $A,B$ have disjoint supports, as do $C,D$, so the coefficient of the original tableau is $1$.

We claim that  
$$\sum_\ell [A E_\ell^{(1)}][ E_\ell^{(2)}D]=0$$
in $\Delta_{2^b1^{a-b}}(\Z^n)$.
Indeed in  each summand we have to distribute the elements of $E$ among the factors $(\ext^2\Z^n)^{\ot b}\ot (\Z^n)^{\ot(a-b)}$. Since the number of these factors is $a$ and $|E|>a$, there is at least one factor of the form $(\ext^2\Z^n)$ containing two elements of $E$. If we exchange the order of these elements the  sum changes sign. However $\sum_\ell [A E_\ell^{(1)}][ E_\ell^{(2)}D]$ is symmetric in   the variables belonging to
E, which yields our claim over $\Z$ by antisymmetrization and torsion-freeness, and hence our proposition. \end{proof}
\begin{example} Let us  illustrate the recursion implied by our identity in a simple case.
Take 
 \[\begin{array}{|c|c|c|}
\hline
2&3&4\\ \hline
2&2&5\\ \hline \end{array}\]
    If we take $E=(2\  \ 2\ 3\ 4)$, we get the relation 
      \[\ 2\  
\begin{array}{|c|c|c|}
\hline
2&3&4\\ \hline
2&2&5\\ \hline
\end{array}
+
\begin{array}{|c|c|c|}
\hline
2&2&4\\ \hline
2&3&5\\ \hline
\end{array}
+
\begin{array}{|c|c|c|}
\hline
2&2&3\\ \hline
2&4&5\\ \hline
\end{array}
=0, \]
in which the original tableau appears with coefficient 2.

If, however, we take $E=(2\ 2 \ \ 2\ 3\ 4)$, we get the relation  
 \[
\begin{array}{|c|c|c|}
\hline
2&3&4\\ \hline
2&2&5\\ \hline
\end{array}
=-
\begin{array}{|c|c|c|}
\hline
2&2&4\\ \hline
2&3&5\\ \hline
\end{array}
-
\begin{array}{|c|c|c|}
\hline
2&2&3\\ \hline
2&4&5\\ \hline
\end{array}
-
\begin{array}{|c|c|c|}
\hline
2&2&2\\ \hline
3&4&5\\ \hline
\end{array}
 .\]
Applying the same relation to the two intermediate tableaux with $E=(2224)$ and $E=(2223)$, respectively, we obtain 
\[
\begin{array}{|c|c|c|}
\hline
2&2&4\\ \hline
2&3&5\\ \hline
\end{array}+
\begin{array}{|c|c|c|}
\hline
2&2&3\\ \hline
2&4&5\\ \hline
\end{array}=-2\ \begin{array}{|c|c|c|}
\hline
2&2&2\\ \hline
3&4&5\\ \hline
\end{array}
.\]
and so 
\[\begin{array}{|c|c|c|}
\hline
2&3&4\\ \hline
2&2&5\\ \hline \end{array}
=\begin{array}{|c|c|c|}
\hline
2&2&2\\ \hline
3&4&5\\ \hline \end{array}.\]

\end{example}

\section{A direct good-filtration argument}
\label{sec:direct-good-filtration}\subsection{Dominance  order}

For partitions $\lambda\vdash d$ of a fixed integer $d$, write $\mu\succeq\lambda$
if we have $\sum_{i=1}^{a}\mu_i\geq\sum_{i=1}^{a}\lambda_i$
for all $a$. \smallskip

 This is the {\em dominance} partial order.  In the
\emph{reverse dominance order} $\mu\preceq\lambda$ means
$\mu \succeq\lambda$.  Dominance implies
$\ell(\mu)\leq\ell(\lambda)$ whenever $\mu \succeq\lambda$.
 \smallskip

Notice that in this partial order the partition formed by the unique part $d$ is maximum while that formed by $d$ parts equal to 1 is minimum.

\medskip
 
  \subsubsection{Good filtrations} In positive  characteristic rational representations of  $GL(V)$  do not generally decompose into direct sums of irreducibles  but  there is a theory of {\em good filtrations} which can be used to  make computations of invariants.
  
  Recall the usual terminology of good filtrations for rational
$\GL(V)$-modules: 
\begin{definition}\label{goodfi}
A finite-dimensional rational $G=\GL(V)$-module has a good or  $\nabla$  filtration
if it has a finite filtration whose successive quotients are induced
modules $\nabla_\lambda(V)$, including determinant twists for rational dominant weights.

Good    $\Delta$  filtrations can be defined in a similar way.
\end{definition} 
It is also known that, if $V$ has a good $\nabla$ filtration then it has a canonical filtration indexed by decreasing dominance order of partitions, such that  successive quotients are direct sums of modules   $\nabla_\lambda(V)$. A similar statement holds for   $\Delta$  filtration but then the order is reversed.\smallskip

We shall use the
standard facts that direct sums of modules with good filtrations again have
good filtrations and that tensor products of modules with good filtrations
have good filtrations. \smallskip

For $\GL(V)$ the tensor-product Theorem goes back to
Wang \cite{Wang}; see also Donkin \cite{Donkin}, Mathieu \cite{Mathieu}, and
Jantzen \cite[II.4]{Jantzen} for geometric proofs, and  to \cite{depr0} for a combinatorial proof.\medskip

 \smallskip

The notion of good filtration and the Tensor product Theorem extend to the integral case.  

By considering  
  the tensor product Theorem applied to $\bigotimes_j\bigwedge^{\lambda_j}\Z^n$ 
the module  $\nabla_\lambda(\Z^n)$ is the  last quotient of  $\bigotimes_j\bigwedge^{\lambda_j}\Z^n$ for the canonical filtration. \smallskip

The highest weight vector is the {\em canonical tableau $C_\lambda$}, the image of  $\bigotimes_je_1\wedge e_2 \wedge\cdots \wedge e_{\lambda_j}$ with all $i$ in the $i^{th}$ column, the character of $C_\lambda$ is $\prod_ix_i^{\lambda'_i}$ with $\lambda'$ the transpose partition.\smallskip

 \subsection{The exterior algebra}

We start from the basic observation which is purely multilinear.

\begin{lemma}\label{lem:exterior-direct-sum}
For finite-dimensional vector spaces $A,B$ (or free modules) and every $h\ge0$ there is a
natural isomorphism
\[
\bigwedge^h(A\oplus B)
\simeq
\bigoplus_{j=0}^h
\bigwedge^j A\otimes \bigwedge^{h-j}B.
\]
\end{lemma}

\begin{proposition}\label{prop:wedge-copies-good}
For every $m,h\ge0$, the $\GL(V)$-modules
\[
\bigwedge^h(V^{\oplus m})
\qquad\text{and}\qquad
\bigwedge^h((V^*)^{\oplus m})
\]
have good filtrations.
\end{proposition}

\begin{proof}
We argue by induction on $m$. The case $m=0$ is immediate, since only the degree-zero component is nonzero. For $m=1$, each nonzero exterior power
$\bigwedge^jV$ is a fundamental highest-weight module; in fact it is an
induced module of minuscule highest weight, and hence has a good filtration
of length one. The same is true of $\bigwedge^jV^*\simeq \bigwedge^{\dim V-j}V\otimes \det^{-1}$, with the corresponding
dominant rational highest weight.\smallskip

Assume the assertion known for $m-1$. Applying
Lemma~\ref{lem:exterior-direct-sum} to
$V^{\oplus m}=V\oplus V^{\oplus(m-1)}$ gives
\[
\bigwedge^h(V^{\oplus m})
\simeq
\bigoplus_{j=0}^h
\bigwedge^jV\otimes
\bigwedge^{h-j}(V^{\oplus(m-1)}).
\]
By induction the second factor in every summand has a good filtration, and
so does $\bigwedge^jV$. \smallskip

The tensor-product Theorem therefore gives a good
filtration on every tensor product occurring on the right, and their direct
sum again has a good filtration. This proves the assertion for $V^{\oplus
m}$. Replacing $V$ by $V^*$ gives exactly the same proof for
$(V^*)^{\oplus m}$.
\end{proof}

\begin{corollary}\label{cor:mixed-exterior-good}
Let $U,W$ be finite-dimensional vector spaces. Every homogeneous component
of the algebra in formula \eqref{IE1},
\[
E(U,V,W)=
\bigwedge^*\!\left(W\otimes V^*\oplus V\otimes U^*\right)
\]
has a good filtration as a $\GL(V)$-module. Consequently $E(U,V,W)$, viewed as
a finite-dimensional rational $\GL(V)$-module, has a good filtration.
\end{corollary}

\begin{proof}
Choose bases of $W$ and $U^*$. As $\GL(V)$-modules,
\[
W\otimes V^*\simeq (V^*)^{\oplus \dim W},
\qquad
V\otimes U^*\simeq V^{\oplus \dim U}.
\]
Proposition~\ref{prop:wedge-copies-good} therefore gives good filtrations on
each homogeneous component of
$\bigwedge(W\otimes V^*)$ and $\bigwedge(V\otimes U^*)$. Since
\[
E(U,V,W)
\simeq
\bigwedge(W\otimes V^*)\otimes
\bigwedge(V\otimes U^*),
\]
the tensor-product Theorem gives the result degree by degree, and then for
the direct sum of all degrees.
\end{proof}

\begin{remark}\label{abw-filtration}
In \cite{ABW}  Akin, Buchsbaum and Weyman developed a characteristic-free theory of standard bitableaux for the exterior algebra $\bigwedge (V\otimes W)$ which produces a filtration of $GL(V)\times GL(W)$-submodules  such that the  associated graded module  is formed  by the modules $\nabla_\lambda(V)\otimes \Delta_{\lambda'}(W)$. This refines  Proposition  \ref{prop:wedge-copies-good}. Furthermore, \cite{GrosshansRotaStein}  develops a more general  standard monomial theory for superalgebras  which could be useful for further generalisations.\end{remark}

 \section{Good filtrations and invariant dimensions}
\label{sec:good-filtration-invariants}

We record the standard consequence needed for invariant theory.  First the module $\nabla_0(V)$  is the trivial 1--dimensional module.

\begin{lemma}\label{lem:invariants-good-filtration}
Let $M$ be a finite-dimensional rational $\GL(V)$-module with a good filtration.
Then
\[
\dim_K M^{\GL(V)}=[M:\nabla_0(V)],
\]
the multiplicity of the trivial induced module in a good filtration of $M$.
\end{lemma}

\begin{proof}
For every dominant weight $\lambda$,
\[
\operatorname{Hom}_{\GL(V)}(K,\nabla_\lambda(V))=0\quad(\lambda\ne0),
\qquad
\operatorname{Hom}_{\GL(V)}(K,\nabla_0(V))=K.
\]
Moreover
$\operatorname{Ext}^1_{\GL(V)}(\Delta(0),\nabla_\lambda(V))=0$, with
$\Delta(0)=K$. Thus taking $\GL(V)$-invariants is exact along a good filtration,
and induction on its length gives the formula. See
Jantzen \cite[II.4.13--II.4.16]{Jantzen} or \cite{depr0}.
\end{proof}

The multiplicities in a good filtration are determined by the formal
character, because
\[
\operatorname{ch}M
=
\sum_\lambda [M:\nabla_\lambda(V)]\operatorname{ch}\nabla_\lambda(V)
\]
and the Weyl characters are linearly independent.

\begin{proposition}
\label{thm:field-independence-GL}
Let $K$ be a  field. For every exterior degree $d$,
\[
\dim_K E(U,V,W)_d^{\GL(V)}
\]
depends only on $\dim U,\dim V,\dim W$ and $d$, and not on $K$.
Consequently the graded dimension, and hence the total dimension, of
$E(U,V,W)^{\GL(V)}$ is independent of the   ground field.
\end{proposition}

\begin{proof}
By Corollary~\ref{cor:mixed-exterior-good}, every graded piece has a good
filtration. If $x_1,\ldots,x_n$ denote the torus weights of $V$, then
\[
\operatorname{ch}E(U,V,W)
=
\prod_{j=1}^n(1+x_j^{-1})^{\dim W}
\prod_{j=1}^n(1+x_j)^{\dim U}.
\]
The same formula, with the exterior degree recorded if desired by an
additional variable, is integral and independent of the characteristic.
Hence the multiplicity of $\nabla_0(V)$ in each graded piece is independent of
$K$. Lemma~\ref{lem:invariants-good-filtration} gives the assertion.
\end{proof}

\begin{remark}
The assumption that $K$ is infinite is needed only when one identifies
invariants for the algebraic group $\GL(V)$ with invariants under the group
of $K$-rational points, see page \pageref{pogf}. 

The good-filtration statement itself is a statement
about rational representations of the algebraic group $\GL(V)$.
\end{remark}

\section{The invariant algebra in positive characteristic}
\label{sec:positive-characteristic-invariants}

We now return to
$
E(U,V,W)
$ in formula \eqref{IE1}
   and ask to what extent the characteristic-zero First Fundamental
Theorem remains valid over a field of positive characteristic.\medskip

We use the notations of \S \ref{bqi},
thus as in page \pageref{ilq} we have
  the map $q$  of formula \eqref{ilq} with image
\[
C_K
=
K[t_{i,j}]
\subseteq
E(U,V,W)^{\GL(V)}
\]
  the subalgebra generated by these quadratic invariants.\smallskip

In characteristic zero, the skew First Fundamental Theorem proved
above gives
\[
C_K
=
E(U,V,W)^{\GL(V)}.
\]

\subsection{A first indication}  Assume
\[
\dim U=\dim W=2, \dim V=4
\]
Then \(C_K\) is generated by the four elements \(t_{i,j},\ i,j=1,2\).
We give \(C_K\) the grading in which
\[
\deg t_{i,j}=1,
\]
so that degree \(d\) below corresponds to exterior degree \(2d\) in
\(E(U,V,W)\).\smallskip

A simple direct computation yields 
the following Hilbert function of the ring generated by $t_{i,j}$
\[
\begin{array}{c|ccccccccc}
\text{degree } d
&0&1&2&3&4&5&6&7&8\\
\hline
\operatorname{char}0
&1&4&10&20&35&20&10&4&1\\
\operatorname{char}2
&1&4&6&4&1&0&0&0&0\\
\operatorname{char}3
&1&4&10&16&19&16&9&0&0\\
\operatorname{char}5
&1&4&10&20&35&20&10&4&1
\end{array}
\]

Thus the algebra generated by the characteristic-zero quadratic
generators already depends on the characteristic.  In this example
characteristic \(5\) happens to give the same Hilbert function as
characteristic zero, whereas characteristics \(2\) and \(3\) do not.\medskip

This is to be contrasted with Proposition~\ref{thm:field-independence-GL}: the dimension of the \emph{full} invariant space
\[
E(U,V,W)^{\GL(V)}
\]
is independent of the   ground field.  Hence what can fail in
positive characteristic is not the stability of the dimensions of the full invariant
space, but rather the assertion that the reductions of the quadratic
elements \(t_{i,j}\) continue to generate it.\smallskip

\begin{remark}
\label{cor:p-nilpotence}
If \(K\) has characteristic \(p>0\), then every positive even-degree
element \(t\) of an exterior algebra over \(K\) satisfies
\[
t^p=0.\quad \text{
In particular,}\quad
t_{i,j}^p=0\quad\text{for all \(i,j\)}\,.
\]
 
\end{remark}
\begin{proof}
Let
\(t=\sum_{\alpha=1}^kc_\alpha m_\alpha\in A\) be a sum of monomials of positive even degree.

Then \(t^p=\sum_{\alpha=1}^kc^p_\alpha m^p_\alpha=0\) since  $m_\alpha^2=0$  for every  monomial of positive degree.
\end{proof}
Thus the algebra generated by the quadratic invariants may become much
smaller after reduction modulo \(p\).  The characteristic-zero FFT is
therefore not expected to survive in positive characteristic with the
same generators.\smallskip

The direct good-filtration argument proved above shows, however, that the full
invariant algebra has a much more rigid behavior.
\medskip

The following simple observation explains how this failure  may be corrected.

\begin{proposition}
\label{lem:factorial-powers}
Let \(A\) be an exterior algebra over a commutative ring and 
\(t=\sum_{\alpha=1}^kc_\alpha m_\alpha\in A\) be a linear combination of monomials $m_\alpha$ of positive even degree.\smallskip

  Then, for every
\(r\ge1\), we have 
$
t^r=r!\,u_r
$
for some \(u_r\in A\).
\end{proposition}

\begin{proof}
 
Since every \(m_\alpha\) has positive exterior degree,
$
m_\alpha^2=0.
$
Moreover, because the \(m_\alpha\) have even degree, they commute:
\[
m_\alpha m_\beta=m_\beta m_\alpha.
\]
Expanding \(t^r\), every term containing a repeated \(m_\alpha\)
vanishes, while every product of \(r\) distinct monomials occurs once
for each permutation of its \(r\) factors.  Hence
\begin{equation}\label{sedi}
t^r=r!
\sum_{\alpha_1<\cdots<\alpha_r}
c_{\alpha_1}\cdots c_{\alpha_r}
m_{\alpha_1}\cdots m_{\alpha_r}=r!u_r;\  u_r:=t^{[r]}.
\end{equation}
The sum is the required element \(u_r\).
\end{proof}
\begin{definition}\label{leld}
The element \(u_r\)  will be denoted by  $t^{[r]}$  and thought of as a {\em divided power}.\smallskip

We also denote  $Exp(t)=\sum_{i=0}^\infty  t  ^{[i]}$ and observe that this infinite sum is in fact truncated  at $k$ the number of monomial terms of $t$.
\end{definition} In fact the map $t\mapsto t^{[r]}$ behaves as an  ordinary divided power since 
 \begin{equation}\label{exp}
(t+s)^{[r]}=\sum_{i=0}^r  t ^{[i]} s ^{[r-i]},\quad  Exp (t+s) =Exp (t )Exp (s)\,.
\end{equation} Thus the even part of an exterior algebra  is a finite-dimensional quotient  of an ordinary algebra of divided powers, but its generators are not all in degree 1.

\begin{remark}

This is the key point: the subalgebra generated by the quadratic
elements can depend strongly on the characteristic, whereas the
dimension of the full invariant space does not.\smallskip

Thus the characteristic-free statement of Proposition \ref{thm:field-independence-GL} is a statement about
the full invariant space and its graded dimension, not about the persistence of
the characteristic-zero generating set. The proof uses only the elementary
exterior-power decomposition and general good-filtration theory.
\end{remark}

\section{Standard monomials over $\mathbb Z$ and the divided-power dual}

\label{sec:integral-smt-divided-power}\subsection{Dominance order and the equivariant standard-basis filtration}

\subsubsection{Standard monomials\label{stmon}}\quad Let $U $ and  $W$ be   finite-rank free abelian groups. 

Put $a=\operatorname{rank} U$, $b=\operatorname{rank} W$,   and  (after choosing bases of the free abelian groups)    consider the
polynomial coordinate ring of the $a\times b$ matrix $\mathtt U=(u_{i,j})$ of indeterminates $u_{i,j}$,
\begin{equation}\label{azz}
A_{\mathbb Z}=\mathbb Z[u_{i,j}]:1\leq i\leq a,\ 1\leq j\leq b 
 =\operatorname{Sym}_{\mathbb Z}(U \otimes W ^{*})=\operatorname{Sym}_{\mathbb Z}(\hom(W,U))\,.
\end{equation}

 We consider the action of the group $G=GL(U)\times GL(W)$  on $A_{\mathbb Z}$ given by the matrix multiplications $g(u_{i,j})h^{-1},\  g\in GL(U),\ h\in  GL(W)$.\ \footnote{this is an example in which $G$ should be thought of as a group scheme since the integer points  are not Zariski dense in the points over $\Q$.}\medskip
 
 Recall the
classical integral standard monomial Theorem for matrices; we use it
without proof.  References include 
\cite{Clausen}, \cite{GrosshansRotaStein},
and \cite{depr0}.

For strictly increasing sequences
$I=(i_1<\cdots<i_k),\ i_k\leq a$ and $J=(j_1<\cdots<j_k),\ j_k\leq b$ of the same length $k\leq \min(a,b)$,
write
\[
 [I\mid J]=\det(u_{i_p,j_q})_{1\leq p,q\leq k}.
\]
A bitableau is a product $\mathtt T=\prod_{r=1}^{s}[I_r\mid J_r]$,
with weakly decreasing row lengths
$\lambda_1\geq\cdots\geq\lambda_s>0$. \smallskip

The {\em shape} of $\mathtt T$ is $\lambda=(\lambda_1,\ldots,\lambda_s)$.\footnote{Sometimes this shape is referred to as {\em row--shape}  and then one has a column shape associated to the dual partition.}  

A bitableau $\mathtt T$ is usually displayed by filling two Young diagrams with the indices $i,j$.  

The polynomial degree of $\mathtt T$    is the number of boxes of each of the two Young diagrams, also called its {\em size}.
\begin{remark}\label{spb}
A basic fact is that the $\Z$ -span of all bitableaux of a given shape  $\lambda$  is a $G$-submodule.
\end{remark}
The bitableau $\mathtt T$  is \emph{standard} if the
row-index tableau and the column-index tableau are semistandard:
entries strictly increase along each row and weakly increase down
columns.

\begin{theorem} 
\label{thm:GRS-standard-  bitableaux }
i)\quad In every degree $d$, the standard   bitableaux  of total size $d$
form a $\mathbb Z$-basis of $(A_{\mathbb Z})_d$.\smallskip

ii)\quad Moreover a bitableau $\mathtt T$ of shape  $\lambda$  is a linear combination of standard bitableaux of shapes $\mu\succeq\lambda$ in the  dominance order.
\end{theorem} This is given by a {\em straightening algorithm} based on certain quadratic relations among determinants of minors (see \cite[Proposition~8.24 and Theorem~8.25]{depr0};
compare \cite[Chapter~11]{BrunsVetter}, \cite{Clausen} and \cite[Chapter~3]{GrosshansRotaStein})
and the equivariant standard filtration
\cite[Definition~9.18, Theorem~9.19 and Corollary~9.20]{depr0}. \medskip

In particular, all the straightening coefficients are integers, and
this basis remains a basis after extension of scalars to any field.\smallskip

We fix an ordering of standard bitableaux, first by shape and then
by any fixed total ordering within each shape.\medskip

Straightening a bitableau of
shape $\lambda$ involves only standard bitableaux of shapes
$\mu \succeq\lambda$, with integral coefficients. The integral
version follows from the integer-coefficient Pl\"ucker straightening
identities and the integral standard-bitableaux bases.\smallskip

\subsubsection{The standard filtration} This Theorem has a representation-theoretic interpretation in terms of good filtrations.\smallskip

Let $\mathcal P_d$ be the finite poset of shapes of standard
bitableaux of degree $d$ for the $a\times b$ matrix $\mathtt U$ ordered by dominance.  If
$\mathcal J\subseteq\mathcal P_d$ is an upper ideal in the
dominance order, that is if  $\lambda\in \mathcal J$ and $\mu\succ \lambda$ then $\mu\in \mathcal J$; set
\begin{equation}\label{eq:ascending-good-filtration}
 F_{\mathcal J,d}:=
 \bigoplus_{\substack{T\text{ standard},\ |T|=d\\
                  \operatorname{shape}(T)\in\mathcal J}}
 \mathbb ZT\ \subseteq A_d.
\end{equation}
The equivariant straightening law makes $F_{\mathcal J,d}$
a $G=\GL(U)\times\GL(W)$-submodule:
the action on minors is expressed by sums of minors of the
same sizes, and standard straightening produces shapes in
$\mathcal J$.  All these modules are saturated $\mathbb Z$-submodules (that is direct summands as groups),
since they are spans of subsets of an integral basis.
No total ordering of shapes is required. \smallskip

 For a single shape
$\lambda \mid |\lambda|=d$ set
\begin{equation}\label{glifl}
F_{\lambda}:=
 \bigoplus_{\substack{T\text{ standard},\ |T|=d\\
                  \operatorname{shape}(T) \succeq \lambda}}
 \mathbb ZT;\qquad F^+_{\lambda}:=
 \bigoplus_{\substack{T\text{ standard},\ |T|=d\\
                  \operatorname{shape}(T)\succ \lambda}}
 \mathbb ZT. 
\end{equation} 
\begin{proposition}\label{grmod}
 The quotient $ F_\lambda/ F^+_{\lambda}$  is a free abelian group with basis the classes of the standard bitableaux of shape $\lambda$, and is a $G$-module.
 
 \begin{equation}\label{clam}
 C_\lambda:=F_\lambda/ F^+_{\lambda}\simeq\nabla_{\lambda}(U)\otimes\nabla_{\lambda}(W^{*})\,.
\end{equation}
\end{proposition}
In particular, when $\lambda=d$ is maximum, $F^+_{d}=\{0\}$  and $C_d$ is the span, in the polynomial algebra $A_{\mathbb Z}$, of the determinants of the $d\times d$ minors of $\mathtt U$.
 \smallskip

\subsubsection{The divided-power algebra}
Now take the \emph{graded $\mathbb Z$-dual} of $A_{\mathbb Z}$ of formula \eqref{azz}
\begin{equation}\label{podi}
 D_{\mathbb Z}:=A_{\mathbb Z}^{\circ}
 :=\bigoplus_{d\geq0}\operatorname{Hom}_{\mathbb Z}
 ((A_{\mathbb Z})_d,\mathbb Z).
\end{equation}
With the multiplication dual to the usual polynomial coproduct
$\Delta(u_{i,j})=u_{i,j}\otimes1+1\otimes u_{i,j}$,
this is the {\em divided-power algebra} on variables $z_{i,j}$.  By duality it carries the contragredient action of $\GL(U)\times\GL(W)$.\smallskip

For a multi-index $\alpha=(\alpha_{i,j})$ set
\[
  u^\alpha=\prod_{i,j} u_{i,j}^{\alpha_{i,j}},
 \qquad
 z^{[\alpha]}=\prod_{i,j}z_{i,j}^{[\alpha_{i,j}]}.
\]
The divided-power monomials form a $\mathbb Z$-basis of $D_{\mathbb Z}$,
characterized by the perfect pairing
\begin{equation}\label{eq:integral-divided-pairing}
 \langle z^{[\alpha]}, u^\beta\rangle=\delta_{\alpha,\beta}.
\end{equation}
In particular,
$z_{i,j}^{[a]}z_{i,j}^{[b]}=\binom{a+b}{a}z_{i,j}^{[a+b]}$;
there is no need to invert any factorial over $\mathbb Z$.\smallskip

For each standard   bitableau   $T$ define $T^\vee\in D_{\mathbb Z}$
to be its element of the \emph{dual standard basis}:
\begin{equation}\label{vee}
 \langle T^\vee,T'\rangle=\delta_{T,T'}
 \qquad\text{for all standard   bitableaux  }T'.
\end{equation}

This dual basis is generally different from a product of formal
determinants in the divided-power variables. The dual standard basis is not in general the divided-power monomial basis. Only its defining
pairing will be used below.\footnote{Note that $D_{\mathbb Z}$ is not a unique factorization domain.}
\begin{remark}\label{ovq} The algebra  $D_\Q=D_\Z\otimes_\Z\Q$ is the usual polynomial algebra in the variables $z_{i,j}$  and then $D_\Z\subset D_\Q$ with $z_{i,j}^{[a]}= \frac{z_{i,j}^ a}{a!}. $

\end{remark}\medskip

\subsection{The decreasing filtration on divided powers}
The partial-order filtration~\eqref{eq:ascending-good-filtration}
induces an order-reversing family of $G$-submodules
\begin{equation}\label{eq:decreasing-dual-filtration}
 F^{\mathcal J}D_d:=
 \{\varphi\in D_d:\langle\varphi,F_{\mathcal J,d}\rangle=0\}
 =\bigoplus_{\substack{T\text{ standard},\ |T|=d\\
              \operatorname{shape}(T)\notin\mathcal J}}
 \mathbb ZT^\vee,
\end{equation}
where $D_d=\operatorname{Hom}_{\mathbb Z}(A_d,\mathbb Z)$.
Thus $\mathcal J\subseteq\mathcal K$ implies
$F^{\mathcal J}D_d\supseteq F^{\mathcal K}D_d$.
For upper ideals $\mathcal J\subseteq\mathcal K$ the perfect
integral pairing yields
\begin{equation}\label{eq:dual-shape-layer}
 F^{\mathcal J}D_d/F^{\mathcal K}D_d
 \simeq(F_{\mathcal K,d}/F_{\mathcal J,d})^*.
\end{equation}
In particular, for $\mathcal J_{\succ \lambda}\subseteq
\mathcal J_{\succeq\lambda}$, by formulas \eqref{glifl} and \eqref{eq:dual-shape-layer}, this quotient is 
\begin{equation}\label{lacl}
\tilde C_\lambda\simeq\nabla_{\lambda }(U )^*\otimes\nabla _{\lambda }(W^* )^*\steq{clam}=\Delta_{\lambda }(U^* ) \otimes\Delta _{\lambda }(W  )
\end{equation}   
\medskip

The filtration is indexed by upper ideals, not by an arbitrarily
chosen total ordering.  \smallskip

When an ordinary finite chain is needed,
a linear extension of the partial order can be taken without
changing any of these submodules or their shape sections.
No assertion is made that $T^\vee$ itself is a product of determinants in
divided powers.\smallskip

 \subsubsection{The generating series}
Over $\mathbb Q$ we have $D_{\mathbb Q}=D_{\mathbb Z}\otimes \Q=\Q[z_{i,j}]$ by Remark \ref{ovq}.

 From here on, we use the transposed contraction convention $t_{i,j}=\sum_r y_{j,r}x_{r,i}$, where $i$ indexes $U$ and $j$ indexes $W$. The skew invariant map  of formula \eqref{ilq}   is, using Definition \ref{leld}
\begin{equation}\label{eq:q-divided-Q}
 q^D :D_{\mathbb Q}\longrightarrow
 E_{\mathbb Q}(U,V,W)^{\mathrm{GL}(V)},
 \qquad q^D :z_{i,j}^{[s]}\longmapsto
 t_{i,j}^{[s]}\steq{sedi}=(\sum_{r=1}^ny_{j,r}x_{r,i})^{[s]}\,,
\end{equation}
   Define
\begin{equation}\label{eq:DivZ-definition}
 \Div_{\mathbb Z}:=q^D (D_{\mathbb Z}).
\end{equation} By formula \eqref{sedi} we have $$q^D (z_{i,j}^{[s]}) =\begin{cases}
\sum_{r_1<\cdots<r_s}y_{j,r_1}x_{r_1,i}y_{j,r_2}x_{r_2,i}\cdots y_{j,r_s}x_{r_s,i}\quad s\leq n\\0\quad s>n
\end{cases}.$$
\begin{equation}\label{ors}\implies\qquad 
\sum_{a=0}^nu_{i,j}^aq^D (z_{i,j}^{[a]}) =\prod_{ r=1}^n(1+u_{i,j} y_{j,r}x_{r,i}),\quad q^D (z_{i,j}^{[n+i]})=0,\ \forall i\geq 1.
\end{equation}
In particular $q^D$  maps $D_\Z$ in  $E_{\mathbb Z}(U,V,W)$ so that $ \Div_{\mathbb Z}\subset E_{\mathbb Z}(U,V,W)^{\mathrm{GL}(V)}.$\medskip

Define the $GL(V)$--invariant exterior-valued forms
\begin{equation}\label{lafd}
\Omega(u) =\sum_{i,j}u_{i,j}\sum_{r=1}^ny_{j,r}x_{r,i}=\sum_{i,j}u_{i,j}q^{ D}(z_{i,j})\,.
\end{equation}
  
All the summands are even and square to zero so we may apply Definition \ref{leld} and formula \eqref{exp}  where we notice that the infinite sums  are truncated at the number of summands of the respective terms. Denote $E_\Z:=E_{\mathbb Z}(U,V,W)$.
Thus  we have  the finite product
\begin{equation}\label{eq:qdefinition}
 \mathcal G (u):=q^{ D}Exp( \sum_{i, j}u_{i,j} z_{i,j})=
 \sum_{\alpha}u^\alpha q^D(z^{[\alpha]})
 =Exp( \Omega (u))= \prod_{r=1}^n\prod_{i=1}^a\prod_{j=1}^b(1+u_{i,j}y_{j,r}x_{r,i})\in A_\Z \otimes E_\Z ,
\end{equation}   where $\alpha$  runs over lists $0\leq \alpha_{i,j}\leq n$ and $u^\alpha=\prod_{i,j}u_{i,j}^{\alpha_{i,j}}$ and $z^{[\alpha]}=\prod_{i,j}z_{i,j}^{[\alpha_{i,j}]}$. \medskip

 This Formula  \eqref{eq:qdefinition} gives
\begin{proposition}\label{intf}
The homomorphism $ q^D :D_{\mathbb Q}\longrightarrow
 E_{\mathbb Q}(U,V,W)^{\mathrm{GL}(V)},$  restricts, for each $d$, to the integral homomorphism $q_d^D:(D_{\mathbb Z})_d\to(E_{\mathbb Z})_{2d}$ 
 specified by the generating identity
\eqref{eq:qdefinition}.
 \end{proposition}  
 
By  \eqref{eq:DivZ-definition}
 $ \Div_{\mathbb Z}:=q^D (D_{\mathbb Z})\subseteq  E_{\mathbb Z}(U,V,W)^{\mathrm{GL}(V)}=E_\Z^{\mathrm{GL}(V)}$.\medskip
 
 Then this homomorphism can be defined over any field of positive characteristic by the same formula and our goal is to prove its surjectivity in Theorem~\ref{thm:Div-base-change-new}.\medskip

A standard shape $\lambda$ is called
\emph{$V$-admissible} when
\begin{equation}\label{eq:admissible-rectangle}
 \ell(\lambda)\leq n,\qquad\lambda_1\leq\min(\dim U,\dim W).
\end{equation}
Thus its diagram is contained in an $n\times m$ rectangle of height $n$  and base $m$, where $m=\min(a,b)$.

These are precisely the shapes needed in the characteristic-zero
invariant decomposition and in the next section.

\section{Integral transpose image and unimodularity}
\label{sec:unimodularity-divided-power}
Throughout this section $A_d=(A_{\mathbb Z})_d$, $D_d=(D_{\mathbb Z})_d=A_d^*$, and $E_{2d}=(E_{\mathbb Z})_{2d}$. 

 \begin{definition}\label{Tabb}
Write  $\mathcal T_{ d}$ for the set of standard bitableaux of degree $d$ and width $\leq \min(a,b)$ and  $\mathcal T_{n,d}\subset \mathcal T_{ d}$ for the set of standard bitableaux with at most $n$ nonempty rows; equivalently their shapes lie in the $n\times\min(a,b)$ rectangle.
\end{definition}    
\par\medskip

\subsection{All coordinates of the transpose}
\label{subsec:all-transpose-coordinates}
 
 We now analyze the transpose $q ^*$ of the integral map $q ^D$  defined by \eqref{eq:qdefinition},   $q ^*:E_{\Z}^*\to A_\Z,\  q_d^*:E_{2d}^*\to A_d$  which is characterised by $\langle q^*(\ell),\varphi\rangle=\langle \ell , q^D\varphi\rangle $.
 
 We now compute its values on the exterior monomial coefficient basis.\smallskip

  If $n=\dim V$ our exterior
algebra has generators given by formula \eqref{levar}
\[
 x_{ r,i}\quad (1\leq i\leq a,\ 1\leq r\leq n),\qquad
 y_{j,r}\quad (1\leq r\leq n,\ 1\leq j\leq b).
\]
Thus $i$ indexes a basis of $U$, $j$ indexes a basis of $W$, and $r$ is the intermediate
index associated with a basis of $V$.\smallskip

 The independent commuting indeterminates
$u_{i,j}$, $1\leq i\leq a$, $1\leq j\leq b$, have no $r$-index.\medskip

Fix the global order for the variables $  x_{ r,i},\ 
 y_{j,r}$ by first increasing intermediate index $r$. 
 
 For each block of a given $r$ we  order first the $y_{j,r}$ in increasing $j$, followed by the $x_{r,i}$ in increasing $i$.\smallskip

 We will write the monomial basis of the exterior algebra in this ordered form and denote it $\mathcal B$.\label{order}\medskip

Consider the finite generating product $
 \mathcal G(u)=\prod_{r=1}^{n}\prod_{i=1}^{a}\prod_{j=1}^{b}
       (1+u_{i,j}y_{j,r}x_{r,i})$ of formula \eqref{eq:qdefinition}.

The factors have even exterior degree and commute; therefore their
ordering does not affect the product. \smallskip

Expanding in the exterior basis $\mathcal B$
produces unique polynomials $F_N(u)\in\mathbb Z[u_{i,j}]$ such that
\begin{equation}\label{tere}
 \mathcal G(u)=\sum_{N\in  \mathcal B} F_N(u)N,\qquad F_N(u)= \langle N^*\mid \mathcal G(u)\rangle.
\end{equation}

\begin{equation}\label{eq:FN-transpose}
F_N(u)=\langle N^*\mid \mathcal G(u)\rangle=\langle N^*\mid q^D(\sum_{\alpha}u^\alpha  
z^{[\alpha]})  \rangle=\langle q^*(N^*)\mid \sum_{\alpha}u^\alpha  
z^{[\alpha]}  \rangle= \sum_{\alpha} u^\alpha  \langle q^*(N^*)\mid   
z^{[\alpha]}  \rangle
\,.\end{equation}

\begin{example}\label{ex:one-index-determinant}
Take $a=b=2$, $n=1$, and abbreviate $x_{1,i }=x_i$ and $y_{ j,1}=y_j$.
The generating product is
\[
 (1+u_{1,1}y_1x_1)(1+u_{1,2}y_2x_1)
 (1+u_{2,1}y_1x_2)(1+u_{2,2}y_2x_2).
\]
For $N=y_1y_2x_1x_2$ only two selections contribute:
\[
 (y_1x_1)(y_2x_2)=-N,\qquad
 (y_2x_1)(y_1x_2)=N.
\]
Consequently
\[
 F_N(u)=-u_{1,1}u_{2,2}+u_{1,2}u_{2,1}
 =-\det\begin{pmatrix}u_{1,1}&u_{1,2}\\u_{2,1}&u_{2,2}\end{pmatrix},
\]
and $q^*(N^*)=-( u_{1,1} u_{2,2}- u_{1,2} u_{2,1})$. \smallskip

The appearance of a determinant will be explained in the next Lemma.
This calculation explains the determinant sign directly from the
anticommutation relations.
\end{example}

\begin{lemma}
\label{lem:allcoords}
For an exterior basis monomial $N$ of bidegree $(d,d)$, and each
$1\leq r\leq n$, let
\[
 I_r=\{i\in\{1,\ldots,a\}:x_{r,i}\text{ occurs in }N\},\qquad
 J_r=\{j\in\{1,\ldots,b\}:y_{j,r}\text{ occurs in }N\}.
\]
If $|I_r|\ne|J_r|$ for some $r$, then $q^*(N^*)=0$.
Otherwise, writing $a_r=|I_r|=|J_r|$, we have
\begin{equation}\label{eq:allcoords}
 q^*(N^*)=\varepsilon_N\prod_{r:\,a_r>0}[I_r\mid J_r],
 \qquad\varepsilon_N\in\{1,-1\}.
\end{equation}
There are at most $n$ nonempty factors. Here $[I_r\mid J_r]$
denotes the determinant of the minor of the matrix  $( u_{i,j})$ with row set $I_r$ and column set $J_r$,
both in increasing order; an empty minor equals $1$.
\end{lemma}
\begin{proof}
Fix an exterior basis monomial $N$ of bidegree $(d,d)$.
For each $r$, write
\[
 I_r=\{i_1<\cdots<i_{s_r}\},\qquad
 J_r=\{j_1<\cdots<j_{t_r}\},
\]
and set
\begin{equation}\label{Nsr}
N_r=
     y_{j_1,r}\cdots y_{j_{t_r},r}x_{r,i_1}\cdots x_{r,i_{s_r} },
\end{equation}
 with empty products equal to $1$.
The product  $N_1\cdot\ldots \cdot N_n$ is well ordered and contains exactly the exterior
variables occurring in $N$. Thus, by the block ordering fixed above, there is no additional
permutation sign:
\begin{equation}\label{eq:N-blocks}
 N=N_1\cdots N_n.
\end{equation}

Recall that
\[
 \mathcal G(u)=\prod_{r=1}^{n}\mathcal G_r(u),\qquad
 \mathcal G_r(u)=\prod_{i=1}^{a}\prod_{j=1}^{b}
 (1+u_{i,j}y_{j,r}x_{r,i}).
\]
Every nonconstant term of each factor of $\mathcal G_r(u)$ contains exactly
one $x$ and one $y$ with intermediate index $r$. 

Therefore, if
$s_r\ne t_r$ for some $r$ in formula \eqref{Nsr}, no term   of the expansion \eqref{tere} of $\mathcal G(u)$ can
contain precisely the exterior variables occurring in $N$, and
the coefficient of $N$ is zero.\smallskip

Now suppose $s_r=t_r=a_r$ for every $r$. Fix $r$ and write
\[
 I_r=\{i_1<\cdots<i_{a_r}\},\qquad
 J_r=\{j_1<\cdots<j_{a_r}\}.
\]
A term  of the expansion \eqref{tere}   contributing to $N_r$ must select $a_r$ nonconstant
factors in $\mathcal G_r(u)$.

By formula \eqref{Nsr} those factors must use each row
index in $I_r$ and each column index in $J_r$ exactly once.
They are consequently indexed by permutations
$\sigma\in S_{a_r}$, with corresponding products
\[
 \prod_{p=1}^{a_r}
 \bigl(u_{i_p,j_{\sigma(p)}}y_{j_{\sigma(p)},r}x_{r,i_p}\bigr).
\]
The selected factors have even exterior degree and hence can
be arranged by increasing $i_p$ without introducing a sign.
Moving all $y$'s to the left of the $x$'s requires
$\binom{a_r}{2}$ transpositions of odd generators; ordering the
$y$'s increasingly requires the further sign
$\operatorname{sgn}(\sigma)$. Thus
\[
 \prod_{p=1}^{a_r}(y_{j_{\sigma(p)},r} x_{r,i_p})
 =(-1)^{\binom{a_r}{2}}\operatorname{sgn}(\sigma)N_r.
\]
Summing over $\sigma$ shows that the coefficient of $N_r$ in
$\mathcal G_r(u)$ is
\(
 (-1)^{\binom{a_r}{2}}
 \det(u_{i,j})_{i\in I_r,\,j\in J_r}
\),
and this also holds for $a_r=0$, using the empty determinant $1$.\smallskip

The exterior generators in different $r$-blocks are disjoint,
so the coefficient of $N_1\cdots N_n$ in $\mathcal G(u)$ is
\[
 (-1)^{\sum_{r=1}^n\binom{a_r}{2}}
 \prod_{r:a_r>0}\det(u_{i,j})_{i\in I_r,\,j\in J_r}
\]
Since $N=N_1\cdots N_n$, setting
\(\varepsilon_N=(-1)^{\sum_{r=1}^n\binom{a_r}{2}}\in\{1,-1\},\)
formula \eqref{eq:FN-transpose} gives
\begin{equation}\label{foco}
q^*(N^*)=\varepsilon_N\prod_{r:a_r>0}[I_r\mid J_r].
\end{equation}
 
There is at most one nonempty factor for each
$r\in\{1,\ldots,n\}$, hence at most $n$ factors.
\end{proof}

\begin{example}\label{ex:two-index-product}
For $a=b=n=2$, use the block-ordered monomial
\[
 N=(y_{1,1}y_{2,1}x_{1,1}x_{1,2})(y_{1,2}x_{2,1}).
\]
The index sets are $I_1=J_1=\{1,2\}$ and
$I_2=J_2=\{1\}$. The first block contributes
$-(u_{1,1}u_{2,2}-u_{1,2}u_{2,1})$, and the second contributes $u_{1,1}$.
Hence, with the indicated block ordering,
\[
 F_N(u)=-u_{1,1}(u_{1,1}u_{2,2}-u_{1,2}u_{2,1}),\qquad
 q^*(N^*)=- u_{1,1}[\{1,2\}\mid\{1,2\}].
\]
Although $u_{1,1}$ can occur in both blocks, there is no extraneous
factor $2$. This illustrates the divided-power normalization.
\end{example}

Recall that $\mathcal T_{n,d}$ is the set of standard bitableaux of
total degree $d$ in the even polynomial coordinates $ u_{i,j}$,
with at most $n=\dim V$ nonempty minor factors (rows of their shapes).

\begin{lemma}\label{lem:distinguished}
For each
$T=\prod_{r=1}^s[I_r\mid J_r]\in\mathcal T_{n,d}$, choose the
exterior monomial
\[
 M_T=\prod_{r=1}^{s}\bigl(y_{ J_r,r}x_{r,I_r}\bigr),\quad x_{r,I_r}=x_{r,i_1 }\cdots x_{r,i_a };\ 
   y_{ J_r,r}=y_{ j_1,r}\cdots y_{ j_a,r}
\]
notice that it respects the prescribed order defined in page \pageref{order}. We then have:\begin{equation}\label{eq:distinguished}
 q^*(M_T^*)=\varepsilon_T T,\qquad
 \varepsilon_T=(-1)^{\sum_{r=1}^s\binom{\lambda_r}{2}}.
\end{equation}
The monomials $M_T$ are pairwise distinct as $T$ varies over standard
bitableaux.
\end{lemma}
\begin{proof}
Apply Lemma~\ref{lem:allcoords} with $I_r,J_r$ equal to the
row and column index sets of the $r^{th}$ minor of $T$. For the fixed global block ordering the additional global-order sign is $1$,
and each block contributes $(-1)^{\binom{\lambda_r}{2}}[I_r\mid J_r]$.
Multiplying gives \eqref{eq:distinguished}. Furthermore, from the
support of $M_T$ one can recover separately the ordered sets
$I_r$ and $J_r$ for each $r$; these determine every minor and hence
$T$. Consequently distinct standard bitableaux yield distinct
exterior basis monomials (up to their ordering signs).
\end{proof}

As a simple consequence we get,  
\begin{theorem}
\label{thm:unimodularity-saturation}\label{thm:split}
Let $m=\min(a,b)$, and write $\mathcal T_{n,d}$ for the standard
bitableaux of degree $d$ whose shapes lie in the $n\times m$ rectangle.
Then
\begin{align}
\label{eq:kernel}
\operatorname{im}q_d^* =\bigoplus_{S\in  \mathcal T_{n,d}}\mathbb Z S;\quad \ker q_d^D&=\bigoplus_{S\in\mathcal T_{ d}\setminus \mathcal T_{n,d}}\mathbb Z S^\vee,\\
\label{eq:Div-basis-admissible}
\Div_{\mathbb Z,d}:=\operatorname{im}q_d^D
 &=\bigoplus_{T\in\mathcal T_{n,d}}\mathbb Z\,q_d^D(T^\vee).
\end{align}
The image is a  direct summand of $E_{2d}$.

\end{theorem}
\begin{proof} By definition, the duals \(N^*\) of the exterior monomials \(N\) of total degree \(2d\) form a basis of \(E_{2d}^*\). By Lemma~\ref{lem:allcoords}, their images in $\Z[ u_{i,j}]$  are up to sign,  products of at most $n$ minors (or zero).
\smallskip 
 
 By integral dominance-triangular straightening
(\cite[Proposition~8.24 and Theorem~8.25]{depr0});  
a product of $s$ minors, of shape
$\lambda$ straightens into an integral linear combination of
standard products of shapes $\mu \succeq\lambda$. Since
$\ell(\mu)\le\ell(\lambda)=s$, none has more than $s$
nonempty factors. Applying this for $s\le n$ gives
$\operatorname{im}q^*\subseteq \bigoplus_{S\in  \mathcal T_{n,d}}\mathbb Z S$. \smallskip
 
 Conversely by Lemma~\ref{lem:distinguished} every $T\in\mathcal T_{n,d}$ lies in  $\operatorname{im}q^*$. 
 
 Thus $\bigoplus_{S\in  \mathcal T_{n,d}}\mathbb Z S\subseteq\operatorname{im}q^*$. 
 Both inclusions are integral and give equality.
 Thus the  image of 
$q_d^*$ is exactly the span of admissible standard bitableaux.\smallskip

The kernel of $\ q_d^D$  is the orthogonal of  the image  of its transpose $q^*$ and this gives
\eqref{eq:kernel}.\smallskip

Since the elements $S^\vee,\ S\in \mathcal T_d$ form a basis  of  $D_d$  formula \eqref{eq:Div-basis-admissible} follows.\medskip

It remains to prove that the image of $q_d^D$ is a direct
summand of $E_{2d}$. Set
\[
B_{n,d}:=\bigoplus_{T\in\mathcal{T}_{n,d}}\mathbb{Z}T.
\]
By formula \eqref{eq:kernel}
\[
q_d^*:E_{2d}^*\longrightarrow B_{n,d}.
\]
is surjective.
Furthermore, by Lemma~\ref{lem:distinguished}, for every $T\in\mathcal{T}_{n,d}$ there is
an exterior basis monomial $M_T$ such that
\[
q_d^*(M_T^*)=\varepsilon_T T,
\qquad \varepsilon_T\in\{1,-1\}.
\]
Thus the assignment $s(T)=\varepsilon_T M_T^*$ extends
linearly to a homomorphism
\[
s:B_{n,d}\longrightarrow E_{2d}^*
\]
with $q_d^*\circ s=\operatorname{id}_{B_{n,d}}$.\smallskip

Dualizing  we deduce that 
$(q_d^*)^*:B_{n,d}^*\longrightarrow E_{2d}$
is a split injection.

Now  $q_d^D:A_d^*=D_d\longrightarrow E_{2d}$
factors as
\[
A_d^*\xrightarrow{\ \rho\ }B_{n,d}^*
\xrightarrow{\ (q_d^*)^*\ }E_{2d},
\]
where $\rho$ is the dual of the inclusion
 $B_{n,d}\subset A_d$. Since  $B_{n,d}$ is spanned by a subset of the
standard-bitableaux basis of $A_d$,  $\rho$ is surjective. Thus,
\[
\operatorname{im}q_d^D=\operatorname{im}(q_d^*)^*,
\]
and everything follows.

 \end{proof}
 For a field $K$ set   $G_K=\GL(U_K)\times\GL(W_K)$. By base change, we obtain the following.
\begin{theorem}
\label{thm:Div-base-change-new} The kernel of the map  $q^D:D_\Z\to \Div_{\mathbb Z}$ is the ideal $I_n$  with $\Z$-basis consisting of the dual basis elements $T^\vee$ whose shapes are not contained in the rectangle $n\times m$, where $m=\min(a,b)$.\smallskip

For every   field $K$, the natural morphism
\begin{equation}\label{eq:Div-base-change-new}
\Div_{\mathbb Z}\otimes_{\mathbb Z}K
\xrightarrow{\ \sim\ } E_K(U,V,W)^{\GL(V_K)}
\end{equation}
is an isomorphism of graded $G_K$-algebras, with divided-power degree $d$ corresponding to exterior degree $2d$.
\end{theorem}
\begin{proof} The first statement is just \eqref{eq:kernel} in Theorem  \ref{thm:unimodularity-saturation}.\smallskip

Integral splitting gives injectivity after tensoring with $K$.
The image consists of invariants. Its graded dimension is independent
of $K$, by the integral standard-monomial basis, and agrees with the
full invariant dimension by Proposition~\ref{thm:field-independence-GL}.
Thus the map is surjective in each homogeneous degree.
\end{proof}
\begin{remark}\label{defn}
The ideals $I_n$ form a decreasing filtration of $G$-stable ideals.
\end{remark}
\subsection{Generators and the Weyl filtration}
\begin{theorem}
\label{thm:generators-all-characteristics}
For any   field $K$ and $n=\dim V$, the invariant algebra
$E_K(U,V,W)^{\GL(V_K)}$ is the quotient $D_K/(I_n\otimes K)$.
It has the standard monomial basis indexed by the dual standard
bitableaux of shapes in the $n\times\min(a,b)$ rectangle, and is
generated as a $K$-algebra by
\[
\{t_{i,j}^{[s]}:1\leq i\leq a,\ 1\leq j\leq b,\ 1\leq s\leq n\}.
\]
In particular, if $\operatorname{char}K=0$, or if
$\operatorname{char}K=p>n$, it is generated by the quadratic
contractions $t_{i,j}=t_{i,j}^{[1]}$.
\end{theorem}
\begin{proof}
Base change gives the asserted quotient and basis.
The divided-power algebra $D_K$ is generated by its one-variable
divided powers $z_{i,j}^{[s]}$ for $s\geq1$; their images are
$t_{i,j}^{[s]}$. Formula \eqref{ors}
gives
$t_{i,j}^{[s]}=0$ for $s>n$. This proves the generating assertion.
When $s!$ is invertible for all $1\leq s\leq n$,
$t_{i,j}^{[s]}=t_{i,j}^s/s!$, so the degree-one divided-power
contractions alone generate.
\end{proof}
\begin{theorem}
\label{thm:Weyl-filtration}
For every $d$, the invariant lattice
$\Div_{\mathbb Z,d}\cong D_d/I_{n,d}$ has a filtration by
$\GL(U_{\mathbb Z})\times\GL(W_{\mathbb Z})$-submodules,
with associated graded the direct sum of the Weyl modules
\begin{equation}\label{wwm}
\Delta_{\lambda }(U_{\mathbb Z}^{*})\otimes
\Delta_{\lambda }(W_{\mathbb Z})
\end{equation}
for each partition $\lambda\vdash d$ contained in the
$n\times\min(a,b)$ rectangle.  The action is the contragredient
of the natural action on $A_{\mathbb Z}=\operatorname{Sym}(U_{\mathbb Z}\otimes W_{\mathbb Z}^{*})$.
This integral Weyl filtration commutes with arbitrary base change.
\end{theorem}
\begin{proof} This follows from formula \eqref{lacl}.
\end{proof}

  \label{part:sl-so}

\section{$\SL(V)$-invariants for vectors}
\subsection{The notion of a skew determinant} In order to analyze the exterior invariant theory of the special linear group $SL(V)$  we need to develop the formalism of 
  determinants  for a square $r\times r$    matrix whose entries are skew-commutative. \smallskip
  
  For an ordinary square matrix $A$ with entries in a commutative ring one has the two usual expansions  by rows or by columns.
  $$\det(A)=\sum_{\sigma\in S_r} \e_\sigma a_ {\sigma(1),1} a_ {\sigma(2),2}\cdots   a_ {\sigma(r),r}=  \sum_{\sigma\in S_r} \e_\sigma a_ {1,\sigma(1) } a_ {2,\sigma(2) }\cdots   a_ {r,\sigma(r) } $$ where $ \e_\sigma$ denotes the sign of the permutation.
  
  These two formulas  when applied  to  a square $r\times r$    matrix, whose entries are skew-commutative, give rise to two different  elements of the exterior algebra which we may call $\overline{\det}_c(X),\  \overline{\det}_r(X)$.\smallskip

Moreover in this exterior case, notice  that    $\overline{\det}_c(X)$ is symmetric, rather than antisymmetric, under column exchange (and  $\overline{\det}_r(X)$ is symmetric   by row exchange) and this permits repeated columns (or rows).  When we have repeated columns one should take into account the multiplicities.

Take two sequences $I=(i_1<\cdots < i_r)$ and $J=(j_1\leq\cdots \leq j_r).$   For each $h$, we let $m_J(h)=\#\{s\mid j_s=h\}$ and define
\[M_c(X)=\prod_hm_J(h)!\]

In the universal exterior algebra over $\Q$, for the $r\times r$ matrix $X:=(x_{i_s,j_t})$ whose entries are skew-commutative, we define 
\begin{equation}\label{skewdet} { {\det}_c}(X)=(I\mid J)_c=\frac{\sum_{\sigma\in S_r}\varepsilon_\sigma\prod_{\ell=1}^rx_{i_{\sigma(\ell)},j_\ell}}{M_c(X)}=\frac{\sum_{\sigma\in S_r}\prod_{\ell=1}^rx_{i_\ell,j_{\sigma(\ell)}}}{M_c(X)}\end{equation}
The integral expression established below defines this element in every characteristic. We have
\begin{lemma}  ${\det}_c$ is a linear combination with integer coefficients of monomials in the $x_{i_s,j_t}$.

We call this the (column) skew determinant. \end{lemma}
\begin{proof} Assume 
that the sequence 
\((j_1,\ldots ,j_r)\)  contains $s$ distinct indices $h_1\,h_2\ldots ,h_s$, $\sum_{p=1}^s m_J(h_p)=r$. Then by \eqref{skewdet}, we get
\[\det\nolimits_c(X)=\sum_{p=1}^s x_{i_1,h_p}\det\nolimits_c(X_p)\]
where $X_p$ is the matrix obtained from X by deleting the first row and any column of index $h_p$. The claim follows immediately by induction.
\end{proof}

In a completely analogous fashion we can define $\det_r(Y):=\det_c(^tY)$ and in this case allow  repeated rows.
If $A$  is  a     matrix of variables commuting with $x_{i,j}$  and with one another:
\begin{equation}\label{molt}
{\det}_c(A\cdot X)=\det(A) {\det}_c(X),\quad  {\det}_r(X\cdot A)= {\det}_r(X)\det(A)\,.
\end{equation}

If $X,Y$ are  square $r\times r$ matrices as above, then   the entries of  $X\cdot Y=T=(t_{i,j})$   commute.  Notice that $T$  has possibly repeated rows and repeated columns. 
Over $\Q$, setting \[\perm(T)=\frac{\sum_{\sigma\in S_r} t_ {\sigma(1),1} t_ {\sigma(2),2}\cdots t_{\sigma(r),r}}{M_r(X)M_c(Y)},\]
it is easy to deduce that 
\begin{equation}\label{perm}
\perm(T)=(-1)^{\binom{r}{2} }  {\det}_r(X) {\det}_c(Y)
\end{equation} 
Indeed,
expanding the permanent and using the anticommutation
of the entries of $X$ and $Y$, we move all the entries
of $X$ to the left, obtaining the sign
$(-1)^{\binom r2}$. The resulting alternating sums
factor as $\det_r(X)\det_c(Y)$. Furthermore, the multiplicities
of repeated rows and columns are exactly cancelled
by the normalizing factors. Formula \eqref{perm} supplies the integral interpretation of the normalized permanent in every characteristic.

\subsection{The module $\Delta_\lambda (W)$ }
Let $K$ be an infinite field of arbitrary characteristic, let $\dim V=n$ and
$\dim W=m$ and put $E=\bigwedge(V\otimes W)=\bigwedge(V^m)$.  \smallskip

As in  formula \eqref{lafond} fix bases $v_1,\ldots,v_n$ of $V$ and $w_1,\ldots,w_m$ of $W$.  We may then think of the set $V^m$ of $m$--tuples of vectors of $V$ as the columns of an $n\times m$ matrix $X$.

Given $n$ column indices  we denote by $[i_1,i_2,\cdots,i_n]$  the determinant $\det_c$  of the square matrix whose columns are columns $i_j$ of $X$, this includes possible repetitions. Notice that $[i ,i ,\cdots,i ]=\prod_{j=1}^nx_{j,i}.$\smallskip

From formula \eqref{molt} it follows:
\begin{proposition}\label{detinv}
The determinants  $[i_1,i_2,\cdots,i_n]$ are $S L(V)$ invariants.
\end{proposition} 

Furthermore:

\begin{lemma}
For $0\le r\le m$, set
\begin{equation}\label{defpro}
z_r=\prod_{j=1}^{r}\prod_{i=1}^{n}x_{i,j}= \prod_{j=1}^{r}[j,j ,\cdots,j]\in\bigwedge^{nr}(V\otimes W),
\end{equation}
 
with the factors taken in the displayed order.  Then, as a $\GL(W)$-module,
\[
  \Span_K\{g\cdot  z_r:g\in\GL(W)\}\simeq  \Delta_{r^n}(W).
\]
\end{lemma}
\begin{proof}
Consider the multidegree $(r,\ldots,r)$ component relative to
$V\otimes W=\bigoplus_{i=1}^{n}(Kv_i\otimes W)$. The exterior-algebra
factorization gives a canonical $\GL(W)$-equivariant isomorphism
\[
 E_{(r,\ldots,r)}\simeq
 \bigotimes_{i=1}^{n}\bigwedge^r(Kv_i\otimes W)
 \simeq (\bigwedge^r W)^{\otimes n}.
\]
Reordering the exterior factors by $i$ identifies $z_r$, up to a harmless
sign, with
\[
 (w_1\wedge\cdots\wedge w_r)^{\otimes n}.
\]
This is the canonical highest-weight vector of weight $(n^r)$ in
$(\bigwedge^r W)^{\otimes n}$. By the defining (integral) realization
of the Weyl module as the submodule generated by this vector, the
$\GL(W)$-module generated by $z_r$ is exactly
$  \Delta_{r^n}(W)$. 
\end{proof}

\begin{theorem}\label{vetinve}
For every infinite field $K$,
\[
 \quad E^{\SL(V)}
 =K\big[(\bigwedge^n(V\otimes W))^{\SL(V)}\big].\quad
\]
In particular the invariant space in degree $d$ is zero if $n\nmid d$.\smallskip

For $0\le r\le m$,
\[
 (\bigwedge^{nr}(V\otimes W))^{\SL(V)}
 \simeq  \Delta_{r^n}(W),
\]
is spanned by products of $r$ maximal skew determinants which we  can   display by a rectangular tableau of shape $r\times n$. 

The products represented by tableaux whose transposes are semistandard form a $K$-basis of $(\bigwedge^{rn}(V\otimes W))^{\SL(V)}$.
\end{theorem}
\begin{proof}
Reasoning as  in Proposition \ref{prop:wedge-copies-good} we  deduce that, for each \(0\leq d\leq nm\),
the  $\SL(V)$-modules \(\bigwedge^d(V\otimes W)\) have good filtrations.

Consequently the dimension of $\SL(V)$-invariants is determined by the
formal character and is independent of the characteristic. Computing it
in characteristic zero by the skew Cauchy identity gives
\[
 \dim_K(\bigwedge^{nr}(V\otimes W))^{\SL(V)}
 =\dim_K  \Delta_{r^n}(W),
\]
and gives zero in all other degrees (since a trivial $\SL(V)$-constituent corresponds to a rectangular $\GL(V)$-partition).

The vector $z_r$ of the lemma is $\SL(V)$-invariant: each block
$\prod_{i=1}^{n}x_{ij}$ transforms by $\det(g)=1$ under
$g\in\SL(V)$. Therefore
\[
   \Delta_{r^n}(W)\simeq K\GL(W)\cdot z_r
 \subseteq (\bigwedge^{nr}(V\otimes W))^{\SL(V)}.
\]
The dimensions agree, so this inclusion is equality.

Finally set $d_j=\prod_{i=1}^{n}x_{ij}\in E_n^{\SL(V)}$, so that 
  $z_r=d_1\cdots d_r$. The subalgebra generated by $E_n^{\SL(V)}$
is stable under $\GL(W)$ and contains $z_r$, hence it contains the
whole cyclic module $K\GL(W)\cdot z_r=E_{nr}^{\SL(V)}$ for every $r$.
These are all the invariant components, proving the assertion.

The last assertion then  follows from Theorem \ref{stat}.(2).
\end{proof}

\section{$\SL(V)$-invariants for vectors and covectors}
\label{sec:slmixed-11nov}

\color{black}\subsection{Invariants}
 Consider as in page \pageref{lafond}\begin{equation}\label{lafond1}
A:=V^a\oplus V^{*b}=V\otimes U^*\oplus W\otimes V^*=\hom(U,V)\oplus \hom(V,W)\,.
\end{equation}   Choosing a basis for $V$, we  display the $a$ vectors in $V^a$  as the columns $x_i=(x_{1,i}, x_{2,i},\dots, x_{n,i})$ of an $n\times a$ matrix $X$  and the $b$   covectors in $V^{*b}$   as the rows $y_i=(y_{i, 1},  y_{ i,2},\cdots, y_{ i,n})$ of a  $b\times n$ matrix $Y$.\medskip

 Let us denote by  $$M=(\bigwedge A)^{SL(V)}=(\bigwedge(V^a\oplus V^{*b}))^{SL(V)}=\bigwedge  (x_{1,i}, x_{2,i},\dots,x_{n,i}; y_{i, 1},  y_{ i,2}, \dots,y_{ i,n})^{SL(V)}\,.$$\medskip
 
 The one-dimensional torus $T=GL(V)/SL(V)$ acts on $M$ via the determinant, and so we can decompose 
 $$M=\oplus_{k\in\Z}M^k;\quad M^k:=\{m\in M\mid g\cdot m=\det(g)^km,\ \forall g\in GL(V)\}\,.$$
 
 Given an admissible partition, that is, a partition $\lambda$ contained in the rectangle $n\times \min(a,b)$ and $k\geq 0$ let $\lambda(k)$ contained in the rectangle $n\times(\min(a,b)+k)$ be obtained from $\lambda$ by adding $k$  columns of length $n$.
 \begin{theorem}\label{fonn}\strut 
 \begin{enumerate}\item $M$ is generated by the $GL(V)$ invariants $(\bigwedge A)^{GL(V)}$ and by the skew determinants $\det_c$ of the maximal minors of $X$ and $\det_r$ of the maximal minors of   $Y$.
\item The space   $M^k$ has a canonical  filtration by $GL(U)\times GL(W)$-modules indexed by  partitions $\lambda$ contained in the rectangle $n\times \min(a,b)$ subject to $\lambda_1\leq\min(a-k,b)$ for $k\geq0$ and $\lambda_1\leq\min(a,b+k)$ for $k<0$, so that the associated graded space is the direct sum of the modules $\Delta_{\lambda(k)}(U^*)\otimes \Delta_{\lambda}(W)$  if  $k\geq 0$  or $\Delta_{\lambda}(U^*)\otimes \Delta_{\lambda(-k)}(W)$  if  $k< 0$. 
\item  The product of a skew determinant  $\det_r$ of a  maximal minor of $Y$ and of a skew determinant  $\det_c$ of a  maximal minor of $X$  is the $GL(V)$  invariant given by formula  \eqref{perm}.
 \end{enumerate}

\end{theorem}
\begin{proof} (1)  

By Theorem \ref{vetinve}    the skew determinants $\det_c$ of the maximal minors of $X$ and $\det_r$ of the maximal minors of   $Y$ are invariant under $SL(V)$.\medskip

Consider now $M^k,\ k>0$; add $k$ new copies  of $V^*$  and  call these new row vectors (or covectors)  $y_{-1},\cdots,y_{-k};\ y_{-i}=(y_{-i,1},y_{-i,2},\cdots,y_{-i,n})$ (the choice of negative indices is needed to make them precede all other vectors when  describing the associated tableaux).

 Take $ \Omega_k=\prod_{i=-k}^{-1}\prod_{r=1}^ny_{i,r}$  a product of determinants $[i,i,\cdots,i]$ and an $SL(V)$ invariant in $M^{-k}$;  and consider  the map  $$\pi:M^k \to \bar  M^0 ;\  m\mapsto m\cdot  \Omega_k$$  where $\bar M$ refers to  $V^a\oplus V^{*(b+k)}$.

The map $\pi$ is a linear isomorphism  to the space $\bar M( \Omega_k)$ of $GL(V)$    invariants  in the enlarged
exterior algebra of $V^a\oplus V^{*(b+k)}$ divisible by $ \Omega_k$ and of degree exactly $nk$ in the coordinates   of the  vectors  $y_{-1},\cdots,y_{-k}$.\medskip

We use the generating functions introduced in formulas \eqref{ors} and \eqref{eq:qdefinition}.   We apply formula \eqref{eq:qdefinition}  in the two cases  $b$ and $b+k$:
\begin{equation}\label{eq1}
\begin{aligned}
 \mathcal G(u)&=\prod_{r=1}^n\prod_{i=1}^a\prod_{j=1}^b(1+u_{i,j}y_{j,r}x_{r,i}),\\
 \widetilde{\mathcal G}(u)&=\prod_{r=1}^n\prod_{i=1}^a\prod_{j\in\{-k,\ldots,-1,1,\ldots,b\}}(1+u_{i,j}y_{j,r}x_{r,i})\\
 &=\mathcal G(u)\prod_{i,r}^{a,n}\prod_{j=-k}^{-1}(1+u_{i,j}y_{j,r}x_{r,i}).
\end{aligned}
\end{equation}

we have
\[
\widetilde{\mathcal G}(u)
=\mathcal G(u)\,P(u),
\qquad
P(u)=
\prod_{j=-k}^{-1}\prod_{r=1}^n
\prod_{i=1}^a(1+ u_{i,j}y_{ j,r}x_{ r,i}).
\]

Since the $y_{j,r}$ are exterior variables, all
products containing the same $y_{ j,r}$ twice vanish.
Consequently,
\[
\prod_{i=1}^a(1+u_{i,j}y_{ j,r}x_{ r,i})
=
1+y_{ j,r}\left(\sum_{i=1}^a u_{i,j}x_{ r,i}\right) .
\]
Hence
\[
P(u)=
\prod_{j=-k}^{-1}\prod_{r=1}^n
\left(1+
 y_{ j,r}\left(\sum_{i=1}^a u_{i,j}x_{ r,i}\right)\right) .
\]

The coefficient of $ \Omega_k$ is therefore
\begin{equation}\label{eq:coefficientY}
(-1)^{\binom{nk}{2}}
\prod_{j=-k}^{-1}\prod_{r=1}^n
\left(\sum_{i=1}^a u_{i,j}x_{ r,i}\right).
\end{equation}
Indeed, the sign comes from moving the 
 $y_{j,r}$'s to the left.\medskip

Now consider one fixed $j$. For a monomial
\[
u_{i_1,j}^{h_1}\cdots u_{i_t,j}^{h_t},
\qquad h_1+\cdots+h_t=n,
\]
its coefficient in
\[
\prod_{r=1}^n\left(\sum_{i=1}^a u_{i,j}x_{ r,i}\right)
\]
is
\[
\sum_{\substack{(s_1,\ldots,s_n)\\
\text{distinct permutations of }
\{i_1^{h_1},\ldots,i_t^{h_t}\}}}
x_{1,s_1}\cdots x_{n,s_n}.
\]
This is precisely the maximal skew determinant
$ {det}_c$ associated with the indicated multiset.\smallskip

It follows that every coefficient of $ \Omega_k$ in
$\widetilde{\mathcal G}(u)$ is, up to the common
sign $(-1)^{\binom{nk}{2}}$, a linear combination
of products of $k$ maximal skew determinants
and coefficients of $\mathcal G(u)$.

By the generating property of
\eqref{eq:qdefinition}, these coefficients span
$\bar M( \Omega_k)$. Transporting them through the
isomorphism given by multiplication by $ \Omega_k$
we conclude
that the elements  of  $M^k$  are linear combinations of products of invariants  times $k$  determinants of maximal minors of $X$.

The case $k<0$   follows in the same way with vectors and covectors interchanged. This proves (1).\medskip

(2) We use the same construction as in (1).
Assume first that $k>0$ and add $k$ new covectors
$y_{-1},\ldots,y_{-k}$. Multiplication by $\Omega_k$
identifies $M^k$ with the space $\bar M(\Omega_k)$
of $GL(V)$-invariants which have degree $n$ in
each of the new covectors.

By Theorem~\ref{thm:Weyl-filtration}, the enlarged
space $\bar M^0$ has a filtration by
$GL(U)\times GL(W\oplus K^k)$-modules whose
associated graded space is
\[
\bigoplus_\mu
\Delta_\mu(U^*)\otimes
\Delta_\mu(W\oplus K^k),
\]
where $\mu$ runs through the partitions contained
in the rectangle $n\times\min(a,b+k)$.

This filtration induces
a filtration on $\bar M(\Omega_k)$.

Consider the component of multidegree
$(n,\ldots,n)$ in the new variables of
$\Delta_\mu(W\oplus K^k)$.

By the standard tableau description of Weyl
modules, this component vanishes unless
$\mu$ contains $k$ complete columns of length $n$.
If $\mu=\lambda(k)$, removing these columns
identifies the component, as a $GL(W)$-module,
with $\Delta_\lambda(W)$.

It follows that the associated graded space
of $\bar M(\Omega_k)\simeq  M^k $ is
\[
\bigoplus_\lambda
\Delta_{\lambda(k)}(U^*)\otimes
\Delta_\lambda(W),
\]
where $\lambda$ runs through the partitions
satisfying \ 
\(
\ell(\lambda)\leq n,\quad
\lambda_1\leq\min(a-k,b).
\)
\smallskip

For $k<0$ the same argument applies after
interchanging vectors and covectors. This proves (2).

\medskip

(3)  is just  formula  \eqref{perm}.

 \end{proof}
\part{Other classical groups}

\section{Good filtrations for classical groups}\label{sec:orthogonal-symplectic-positive-characteristic}
We now explain how the preceding good-filtration argument applies to the invariant theory of 
the symplectic and orthogonal groups.\smallskip

Let \(K\) be a field (of characteristic different from \(2\) in the orthogonal case) and let \(V\) be a finite-dimensional
\(K\)-vector space.  In the symplectic case we assume that \(V\) is
equipped with a nondegenerate alternating form and denote by
$Sp(V)$ the symplectic group preserving the form.\medskip

In the orthogonal case we equip \(V\) with a nondegenerate symmetric bilinear form. We denote by
$
 \SO(V),\   \mathrm O(V),
$ the special orthogonal group  and the orthogonal group  preserving the form.\medskip

Recall that a pair \(H\subset G\) of reductive
groups is called {\em a Donkin pair} if every rational \(G\)-module with a
good filtration restricts to an \(H\)-module with
a good filtration.

For the classical inclusions considered here recall (Brundan \cite{brundan} Proposition 3.3(i), (iii), Proposition 3.4(i), and Lemma 3.2)

\begin{theorem}
\label{thm:restriction-good-filtration-classical}
Let \(M\) be a rational \(\GL(V)\)-module having a good filtration.

\begin{enumerate}
\item If \(V\) is symplectic, then the restriction of
\(
M\) to 
\(Sp(V)\)
has a good filtration as an \(\Sp(V)\)-module.

\item If \(V\) is orthogonal and \(\operatorname{char}K\neq2\), then
the restriction of
\(
M\) to 
\(SO(V)\) has a good filtration as an \(\SO(V)\)-module.
\end{enumerate}
\end{theorem}

By Proposition \ref{prop:wedge-copies-good}
for every $m,h\ge0$, the module $
\bigwedge^h(V^{\oplus m})$ has a good filtration
 as a \(\GL(V)\)-module.
Hence Theorem~\ref{thm:restriction-good-filtration-classical} gives
immediately:

\begin{corollary}
\label{cor:classical-good-filtration-exterior}
Let $\dim W=m$. The restriction of \(\bigwedge^h(V^{\oplus m})=\bigwedge^h(V\otimes W)\) to \(\Sp(V)\) has a good filtration.

If \(\operatorname{char}K\neq2\), its restriction to \(\SO(V)\) has a
good filtration as well.
\end{corollary}

We can therefore repeat the base-change argument used for
\(\GL(V)\).\smallskip

In the symplectic case choose integral lattices together with a
unimodular alternating form and let
\(
E_{\mathbb Z}
\)
denote the corresponding integral exterior algebra.  \smallskip

In the
orthogonal case it is convenient to work over
\(
R=\mathbb Z[1/2],
\)
with a unimodular symmetric form, so that the distinction between
\(\SO\) and \(\mathrm O\) behaves uniformly under base change.

\begin{theorem}
\label{thm:classical-field-independence}
Let \(K\) be a field.

\begin{enumerate}
\item In the symplectic case,
\[
\dim_K
\bigwedge^h(V\otimes W)^{\Sp(V_K)}
\]
is independent of \(K\).

\item In the orthogonal case, for
\(\operatorname{char}K\neq2\),
\[
\dim_K
\bigwedge^h(V\otimes W)^{\SO(V_K)}
\]
is independent of \(K\).

\item Under the same assumption
\(\operatorname{char}K\neq2\),
\[
\dim_K
\bigwedge^h(V\otimes W)^{\mathrm O(V_K)}
\]
is independent of \(K\).
\end{enumerate}
\end{theorem}

\begin{proof}
For the connected groups $H=\Sp(V),\ \SO(V)$, restriction
has a good filtration by the Donkin-pair Theorem. For such a module
$M$, the dimension of $M^H$ equals its good-filtration multiplicity
$[M:\nabla_H(0)]$. This is determined by the formal character of $M$:
the   characters of the dual Weyl modules form an integral basis of the characters, while the
formal character of $\bigwedge^h(V\otimes W)$ is independent of the
field. Hence the invariant dimension is constant.\smallskip

To justify the assertion for the full orthogonal group, work over
$R=\mathbb Z[1/2]$ with the chosen unimodular form and put
$M_R=\bigwedge^h(V_R\otimes W_R)$ and
$L_R=M_R^{\SO(V_R)}=  M_R\cap M_\Q^{\SO(V_\Q)}$. \smallskip

By construction $L_R$ is a direct summand  of the $R$-module $M_R$  and, since the dimension of the $\SO(V_K)$ invariants is independent of the field $K$ we have $L_R\otimes_R K=M_K^{\SO(V_K)}$ for every field $K$ of characteristic different from two.
\smallskip

The quotient $\mathrm O(V_R)/\SO(V_R)$ is the   group $C_2=\{1,s\}$ which  acts on $L_R$. 

Since $2$ is invertible, the
idempotent $(1-s)/2$ decomposes  $L_R= L_R^+\oplus L_R^-$ with  $L_R^+=M_R^{\mathrm O(V_R)}$.
It follows that also $M_R^{\mathrm O(V_R)}$ is a free direct summand whose rank equals the
characteristic-zero invariant dimension and everything follows as in the $SO(V)$ case.
\end{proof}

\section{Standard tableaux for symmetric and alternating matrices}
\label{sec:classical-standard-tableaux} Our next goal is to give a description of the rings of invariants  as quotients of suitable divided-power algebras.\medskip

We begin by recalling the standard-tableaux and straightening results
of De Concini--Procesi~\cite{DeConciniProcesi1976} for symmetric and
alternating matrices. Their characteristic-free treatment supplies the
integral straightening relations for products of minors of a symmetric
matrix and of Pfaffians of an alternating matrix.\smallskip

 What we use is that
straightening is homogeneous in the polynomial degree, is triangular
in the dominance order on shapes, and cannot increase the number of
nonempty factors (rows in the factor-size tableau convention). \smallskip

We state the conventions and consequences needed in this paper; the
subsequent good- and Weyl-filtration arguments proceed as in
Sections~7--8.  \smallskip

With the notations of \S \ref{sec:char-zero-orth-symp} let $U=\Z^m=\oplus_j\Z u_j$  and $V=\Z^n=\oplus_j\Z v_j$ with the symmetric form $\sum_{i=1}^nx_i^2$ in the orthogonal case while $V=\Z^{2n }=\oplus_{j=1}^n\Z e_j\oplus \oplus_{j=1}^n\Z f_j$ where    
$e_1,\dots,e_n,f_1,\dots,f_n$ is a symplectic basis. Thus $n$ is the orthogonal dimension or half the symplectic dimension.

For the orthogonal group one can take any symmetric unimodular quadratic form, at the level of invariants this does not change the computations, since any two unimodular integral quadratic forms of the same rank become isometric over the integral elements \(\mathcal O_K[1/2]\) for a suitable number field \(K\) (cf. \cite{O'Meara} \S\S~92, 102--104).\medskip

We study the invariant theory  of the integral forms $\left(\bigwedge^*(V \otimes U )\right)^{G}$ with $G$  the symplectic or the orthogonal group.  The group $\GL(U)$ acts on the invariants. We then apply this to  the invariant theory  of   $\left(\bigwedge^*(V_K \otimes U_K )\right)^{G}$ where $K$ is a field, of characteristic $\neq 2$ in the orthogonal case.\smallskip

\subsubsection{Normalisation of the symplectic diagonal generator.}
In the symplectic case, the diagonal quadratic generators admit an integral normalization that must be taken into account in characteristic two. 
\smallskip

Write
$x_{a,i}=e_a\otimes u_i$, $y_{a,i}=f_a\otimes u_i$ so that $\bigwedge^*(V \otimes U )=\bigwedge^*(x_{a,i},\ y_{a,i})$.

As in Theorem \ref{thm:symplectic} consider the matrix $X$ with rows $x_{1,i},\cdots, x_{n,i},y _{1,i},\cdots,  y_{n,i}$ and set:
  \[
T={}^tXJX,\quad\text{with}\quad  t_{i,j}=\sum_{a=1}^n(x_{a,i}y_{a,j}+x_{a,j}y_{a,i})\,.
\]
The matrix $T$ is symmetric and its entries are invariants.

We have  
 $t_{i,i}=2\sum_{a=1}^n x_{a,i}y_{a,i}:=2q_{i,i} $.
 We therefore replace the matrix $T$ of generators
introduced in Section~\ref{sec:char-zero-orth-symp}
with the integral matrix $Q$:
\begin{equation}\label{imq}
 q_{i,i}=\sum_{a=1}^n x_{a,i}y_{a,i},
\,\quad 
 q_{i,j}=\sum_{a=1}^n(x_{a,i}y_{a,j}+x_{a,j}y_{a,i}),\ i\neq j.
\end{equation}

\subsubsection{Symmetric matrices.}
Let $m=\operatorname{rank} U$, and let $X^+=(u_{i,j})_{1\le i,j\le m}$, $u_{i,j}=u_{j,i}$,
be a generic symmetric matrix. 

 For increasing sequences
$I=(i_1<\cdots<i_p)$ and $J=(j_1<\cdots<j_q)$, write
$I\preceq J$ if $p\ge q$ and $i_t\le j_t$ for $1\le t\le q$. The determinant of an $a\times a$  minor with rows $I$ and columns $J$ is denoted 
$$\Bigg [ {I \atop   J}\Bigg]=\Bigg [ {J\atop  I}\Bigg]\qquad \text{by symmetry},$$
where $|I|=|J|=a$. The two increasing index sequences
form the two rows of the minor tableau. \smallskip

We use the following standardness condition for symmetric minors:
\[
 T=\Bigg [ {I_1\atop  J_1}\Bigg]\cdots\Bigg [ {I_s\atop  J_s}\Bigg],\qquad
 I_1\preceq J_1\preceq I_2\preceq J_2\preceq\cdots
 \preceq I_s\preceq J_s,
\]
where factors are ordered so that $a_1\ge\cdots\ge a_s>0$,
$a_t=|I_t|=|J_t|$. Thus the tableau has two rows of length
$a_t$ for each minor of size $a_t$. Its Young-diagram shape is
\[
 \tau^+(\lambda)=(a_1,a_1,a_2,a_2,\ldots),\qquad
 \lambda=(a_1,\ldots,a_s),\quad |\lambda|=d.
\]
\emph{Every column of this tableau diagram has even height.}

Equivalently, its highest weight $2\eta$ is twice a dominant
polynomial weight $\eta$ (see page \pageref{sec:direct-good-filtration}). In degree one the tableau shape
is $(1,1)$ and the highest weight is $2\omega_1$.\smallskip

\subsubsection{Alternating matrices.}
Let $X^-=(u_{i,j})$ be a generic alternating $m\times m$ matrix with
$u_{j,i}=-u_{i,j}$ and $u_{i,i}=0$. For an increasing
sequence $I=(i_1<\cdots<i_{2a})$, let $[I]_{\mathrm{Pf}}$
be the principal Pfaffian of order $2a$ (polynomial degree $a$).
A standard Pfaffian product of degree $d$ is
\[
 T=[I_1]_{\mathrm{Pf}}\cdots[I_s]_{\mathrm{Pf}},\qquad
 I_1\preceq I_2\preceq\cdots\preceq I_s,
\]
where $|I_t|=2a_t$ and $a_1\ge\cdots\ge a_s>0$.
Each Pfaffian contributes one row of length $2a_t$; thus
\[
 \tau^-(\lambda)=(2a_1,2a_2,\ldots),\qquad
 \lambda=(a_1,\ldots,a_s),\quad |\lambda|=d.
\]
\emph{Every row of this tableau diagram has even length.}
 
Equivalently the highest
weights are precisely the nonnegative integral combinations
$\sum_{i=1}^{\lfloor m/2\rfloor}c_i\omega_{2i}$.
In degree one the tableau shape is $(2)$ and the highest
weight is $\omega_2$.  \smallskip

For integral Pfaffian straightening and the two-factor
laws below, see De Concini--Procesi
\cite[Lemmas~6.1--6.2 and Theorem~6.5]{DeConciniProcesi1976}.\smallskip

In either case the first row of the tableau has length at most
$m$; equivalently the highest-weight partition $\rho=\tau'$
has at most $m$ nonzero parts.\smallskip

In either case the standard tableaux of polynomial degree $d$
form a $\Z$--basis of the degree $d$  part of $\Z[u_{i,j}]$.\smallskip

 For a product of  two minor factors with sizes $a\ge b$,
or two Pfaffians of orders $2a\ge2b$, every standard term in
the straightening has at most two nonempty factors, of respective
sizes $s_1,s_2$ (half-orders in the Pfaffian case), and
\[
 s_1+s_2=a+b,\qquad s_1\ge a.
\]
The second factor may be empty. These are the straightening
laws we use from De Concini--Procesi
\cite{DeConciniProcesi1976}; compare in particular
\cite[Theorem~6.5]{DeConciniProcesi1976} for Pfaffians.\smallskip

Therefore a product of
  shape $\lambda$
straightens only into shapes $\mu \succeq\lambda$; in particular
straightening cannot increase the number of factors.\smallskip

For each $d$, the $\Z$--span of standard tableaux of degree $d$ and of  shapes
belonging to a lower ideal in reverse dominance is a
$\GL(U)$-submodule.   We can keep the notations   of formula \eqref{glifl} except that now $T$ denotes a single tableau and not a bitableau.\smallskip

The quotient $ F_\lambda/ F^+_{\lambda}$  is a free abelian group with basis the classes of the standard tableaux of shape $\lambda$, and is a $G$-module.
 
 \begin{equation} \label{clams1}
 C_\lambda:=F_\lambda/F^+_\lambda\simeq\nabla_{\tau^\epsilon(\lambda)}(U),
\end{equation}
where $\lambda$ is the factor-size shape and $\tau^\epsilon(\lambda)$ is the full tableau shape.
 
\section{Classical contractions, integral splitting, and Weyl filtrations}
\label{sec:classical-divided-powers}
 
\subsection{The coefficient--wise contraction maps}
Define as in formula \eqref{podi}\begin{equation}\label{podie}\text{ the divided power algebras}\quad 
 D_{\mathbb Z}^\epsilon:= \bigoplus_{d\geq0}\operatorname{Hom}_{\mathbb Z}
 ((A_{\mathbb Z}^\epsilon)_d,\mathbb Z).
\end{equation}
where $A^+=\Z[u_{i,j}:i\leq j]$ with $u_{j,i}=u_{i,j}$, and $A^-=\Z[u_{i,j}:i<j]$ with $u_{j,i}=-u_{i,j}$ and $u_{i,i}=0$.

The elements of the divided power algebra will be denoted by $z_{i,j}^{[a]}$ dual to $u_{i,j}^{ a }$.\medskip

In the orthogonal case let $V^-$ be the  quadratic module of dimension $n$
  over $\Z$, with orthonormal basis $e_1,\ldots,e_n$. \smallskip

 Let
$x_{ri}$ ($1\leq r\leq n$, $1\leq i\leq m$) be the exterior generators
of $E^-=\bigwedge(V^-\otimes U)$.  \smallskip

In the symplectic case let $V^+$ be a
rank--$2n$ symplectic module with a symplectic basis, and let
$x_{a,i},y_{a,i}$ ($1\leq a\leq n$, $1\leq i\leq m$) be the generators
of $E^+=\bigwedge(V^+\otimes U)$. \smallskip

 Choose the total exterior order by
increasing $r$ in the orthogonal case, and by increasing $a$ with all
$x_{ai}$ before all $y_{ai}$ in each symplectic block.\medskip

The generators for the symplectic case  are given by the integral matrix $Q$ of formula \eqref{imq}:
\begin{equation}\label{imq1}
 q_{i,i}=\sum_{a=1}^n x_{a,i}y_{a,i},
\,\quad 
 q_{i,j}=\sum_{a=1}^n(x_{a,i}y_{a,j}+x_{a,j}y_{a,i}),\ i\neq j.
\end{equation}
The signs correspond to our choice of bases; changing the sign of the
symplectic form changes the global convention, not the conclusions. \smallskip

For the orthogonal group we will also denote the  generators as $q_{i,j}$ rather than $t_{i,j}$ ($i<j$) in order to maintain uniform notation.\medskip

 Formula \eqref{exp}  shows that there is a homomorphism  $q^\epsilon:D^\epsilon\to  E^\epsilon$ with \begin{equation}\label{eq:classical-q-definition}q^\epsilon(z_{i,j}^{[a]})=q_{i,j}^{[a]};\qquad  
 q^\epsilon(z^{[\alpha]})=q^{[\alpha]},
 \qquad z^{[\alpha]}=\prod_I z_I^{[\alpha_I]}\in D^\epsilon,\ q^{[\alpha]}=\prod_I q_I^{[\alpha_I]}\in E^\epsilon.
\end{equation}Define the $G$--invariant exterior-valued forms
\begin{align}
 \Omega^-(u)&=\sum_{i<j}u_{i,j}\sum_{r=1}^nx_{ri}x_{rj}=\sum_{i<j}u_{i,j}q^{ -}(z_{i,j}),
 \label{eq:classical-omega-minus}\\
 \Omega^+(u)&= \sum_i u_{i,i}\sum_{a=1}^nx_{a,i}y_{a,i}
       +\sum_{i<j}u_{i,j}\sum_{a=1}^n(x_{a,i}y_{a,j}+x_{a,j}y_{a,i})=\sum_{i\leq j}u_{i,j}q^{ +}(z_{i,j}) .
 \label{eq:classical-omega-plus}
\end{align}

All the summands are even and square to zero so we may apply Definition \ref{leld} and formula \eqref{exp}  where we notice that the infinite sums  are truncated at the number of summands of the respective terms. \smallskip

Put $\mathcal I_- =\{(i,j):1\leq i<j\leq m\}$ and $\mathcal I_+=\{(i,j):1\leq i\leq j\leq m\}$. Thus as in formula \eqref{eq:qdefinition} we have the finite product
\begin{equation}\label{eq:classical-generating-product}
\begin{aligned}
 \mathcal G^\epsilon(u)&:=q^\epsilon Exp\left(\sum_{(i,j)\in\mathcal I_\epsilon}u_{i,j}z_{i,j}\right)=Exp(\Omega^\epsilon(u))\\
 &=\prod_{\substack{\text{individual square-zero}\\\text{summands }v\text{ of }\Omega^\epsilon(u)}}(1+v)
 =\sum_\alpha u^\alpha q^\epsilon(z^{[\alpha]})\in A^\epsilon\otimes E^\epsilon.
\end{aligned}
\end{equation}  This Formula  
has integral coefficients; in the displayed sum the monomial $u^\alpha$
uses the independent upper--triangular coordinates (in the orthogonal case the diagonal entries are 0).  The definition is
independent of the order of its even factors.  

\begin{remark}
Here a divided power is defined by the indicated coefficient, not by
inverting a factorial.  The ordinary exponential identities over
$\mathbb Q$ show that $q^\epsilon$ respects divided-power
multiplication; since both source and target are torsion-free over
$\Z$, the same identities hold integrally. 

The forms
$\Omega^\epsilon$ are invariant under the corresponding classical
groups, hence every coefficient is invariant (equivalently use the
polarised divided-power coaction).  \smallskip

The maps are equivariant for
$\GL(U)$ and homogeneous with $q_d^\epsilon:D_d^\epsilon\to E_{2d}^\epsilon$.

For an exterior monomial $M$, $(q_d^\epsilon)^*(M^*)$ is the coefficient of \(M\) in \(\mathcal G^\epsilon(u)\) as follows
 by pairing both sides with $z^{[\alpha]}$.

\end{remark}
In analogy with Lemma \ref{lem:allcoords} we have
\begin{lemma}
\label{lem:classical-transpose}
Use the exterior monomial basis in the
block order fixed above. Write $Z^\pm=(u_{i,j})$ for the polynomial coordinate matrices. Write $M^*$ for the dual of a basis monomial $M\in E_{2d}^\epsilon$.  In the orthogonal case,
write
\[
 M=\prod_{r=1}^n x_{r,I_r},\]
if
$I=\{i_1<\cdots<i_a\}$, 
$x_{r,I}=x_{r,i_1}\cdots x_{r,i_a}$,
with empty product $1$.  

Some $I_r$ may be empty. Then
\begin{equation}\label{eq:classical-all-coordinates}
 (q_d^-)^*(M^*)=
 \begin{cases}
 \displaystyle\prod_{r=1}^n\operatorname{Pf}(Z^-_{I_r}),
       & |I_r|\text{ is even for all }r,\\
 0,&\text{otherwise}.
 \end{cases}
\end{equation}
The Pfaffian of the empty matrix is $1$.\medskip

In the symplectic case write
\[
 M=\prod_{r=1}^n x_{r,I_r}y_{r,J_r},\qquad
\]
where if
$I=\{i_1<\cdots<i_a\}$, (resp.  $J=\{j_1<\cdots<j_b\}$), 
$x_{r,I}=x_{r,i_1}\cdots x_{r,i_a}$ (resp. $y_{r,J}=y_{r,j_1}\cdots y_{r,j_b}$,
with empty product $1$.  
Then
\begin{equation}\label{eq:classical-all-coordinates-plus}
 (q_d^+)^*(M^*)=
 \begin{cases}
 \displaystyle(-1)^{\sum_{r=1}^n\binom{|I_r|}{2}}
 \prod_{r=1}^n\det Z^+_{I_r,J_r},
       & |I_r|=|J_r|\text{ for all }r,\\
 0,&\text{otherwise}.
 \end{cases}
\end{equation}
The determinant of the empty matrix is $1$. \smallskip

 In particular the
nonzero coordinates are products of at most $n$ Pfaffians of principal minors,
respectively at most $n$ determinants of symmetric minors, with no additional
multiplicities. 
\end{lemma}
\begin{proof}
We prove the two identities directly, including the coefficient
normalizations.  \smallskip

By the definition of $q_d^\epsilon$ in
\eqref{eq:classical-generating-product}--\eqref{eq:classical-q-definition},
for every exterior basis monomial $M$ of exterior degree $2d$ its coefficient in $\mathcal G^\epsilon(u)$ equals
\[\langle  M^* ,\mathcal G^\epsilon(u)  \rangle= \langle  M^* , q^{ \epsilon}Exp( \sum_{(i,j)\in\mathcal I_\epsilon}u_{i,j} z_{i,j}) \rangle=   
 \sum_{ \alpha }u^\alpha
    \langle(q ^\epsilon)^*(M^*),z^{[\alpha]}\rangle.
\]
Because the monomial basis of $A_d^\epsilon$ is dual to the
divided-power basis of $D_d^\epsilon$, this proves that
\begin{equation}\label{eq:classical-transpose-coefficient-proof}
 (q_d^\epsilon)^*(M^*)=\sum_{|\alpha|=d}u^\alpha
    \langle(q_d^\epsilon)^*(M^*),z^{[\alpha]}\rangle.
\end{equation}
where the polynomial on the right is viewed in $A_d^\epsilon$.

 \medskip
 
\emph{Orthogonal case.}
Factor the generating product by orthogonal basis index:
\[
 \mathcal G^-(u)=\prod_{r=1}^n G^-_r(u),\qquad
 G^-_r(u)=\prod_{1\le i<j\le m}
             (1+u_{i,j}x_{r,i}x_{r,j}).
\]
Each nonconstant selected factor in $G^-_r$ uses exactly two
$x$--generators of the $r$--th block. Hence a prescribed block
monomial $x_{r,I}$ with $|I|$ odd has coefficient zero.\smallskip

Suppose $I=\{i_1<\cdots<i_{2a}\}$. A nonzero contribution
to the coefficient of $x_{r,I}=x_{r,i_1}\cdots x_{r,i_{2a}}$
must select precisely $a$ factors whose unordered index
pairs form a perfect matching of $I$. 

Conversely, every
perfect matching determines a unique collection of such
factors, with pairwise disjoint indices, and therefore
gives a nonzero exterior product.\smallskip

All selected factors have even exterior degree, so
reordering them does not affect the sign. \smallskip

The sign of
an individual matching is exactly the sign of the
permutation needed to order its $2a$ odd variables as
$x_{r,i_1}\cdots x_{r,i_{2a}}$.\smallskip

The defining matching expansion of the Pfaffian is
\[
 \operatorname{Pf}(Z^-_I)=
 \sum_{\mathfrak m}
  \operatorname{sgn}(\mathfrak m)
  \prod_{\{i,j\}\in\mathfrak m,\ i<j}u_{i,j},
\]
where the sum runs over perfect matchings and the
matching sign is computed with respect to the increasing
order of $I$.

Thus we deduce that
the  coefficient  of   \([x_{r,I}] \) in \(G^-_r(u)\)   equals \(\operatorname{Pf}(Z^-_I),
\)
with every matching occurring exactly once and with
coefficient $1$ or $-1$.

The block variables are disjoint and the block ordering
of $M$ is fixed. Multiplying these identities over $r$
proves \eqref{eq:classical-all-coordinates}.

\medskip
\emph{Symplectic case.}
Similarly, write
\[
 \mathcal G^+(u)=\prod_{r=1}^nG^+_r(u),\qquad
 G^+_r(u)=
 \left(\prod_{i=1}^m(1+u_{i,i}x_{r,i}y_{r,i})\right)
 \left(\prod_{1\le i<j\le m}
     (1+u_{i,j}x_{r,i}y_{r,j})
     (1+u_{i,j}x_{r,j}y_{r,i})\right).
\]
This is the factorization of the product already defining
$q^+$: the two terms involving $u_{i,j}$ are separate even,
square-zero exterior factors. Notice also that the diagonal
factor has coefficient $1$, in accordance with the
normalization of the diagonal symplectic invariant
introduced earlier.

Every selected nonconstant factor uses one $x$ and one
$y$. Consequently, the coefficient of
\(
 x_{r,I}y_{r,J}\) in $G^+_r(u)$ is  $0$
 unless \( |I|=|J|.\)

Suppose
$I=\{i_1<\cdots<i_a\}$ and
$J=\{j_1<\cdots<j_a\}$.
Any selection producing exactly $x_{r,I}y_{r,J}$
must use each $x_{r,i_p}$ and each $y_{r,j_q}$
exactly once.\smallskip

Such selections are in bijection with permutations
$\sigma\in S_a$: for each $p$, select the factor
containing
\[
 u_{i_p,j_{\sigma(p)}}x_{r,i_p}y_{r,j_{\sigma(p)}}.
\]
Indeed, for every ordered pair $(i,j)$ there is exactly
one factor containing the term $u_{i,j}x_{r,i}y_{r,j}$:
the diagonal factor if $i=j$, and one of the two
off-diagonal factors if $i\ne j$, where
$u_{j,i}=u_{i,j}$. Thus every permutation gives
exactly one selection, even when $I\cap J$ is nonempty.\smallskip

Since the selected factors have even exterior degree, we may reorder them so that the indices \(i_p\) appear in increasing order, without changing the sign.\smallskip

 Moving all the $y$'s after all the
$x$'s requires $\binom a2$ interchanges of odd
variables, while sorting the $y$'s into increasing
order contributes $\operatorname{sgn}(\sigma)$.

Therefore the coefficient of \(x_{r,I}y_{r,J}\) in \(G^+_r(u)\) equals
\[
\begin{aligned}
 (-1)^{\binom a2}
   \sum_{\sigma\in S_a}\operatorname{sgn}(\sigma)
     \prod_{p=1}^a u_{i_p,j_{\sigma(p)}} =(-1)^{\binom a2}\det Z^+_{I,J}.
\end{aligned}
\]
The formula also holds for $a=0$.\smallskip

The different $r$-blocks involve disjoint exterior
generators, and the block ordering is fixed.
Multiplication therefore proves
\eqref{eq:classical-all-coordinates-plus}
with exactly the stated sign. The restriction on the
number of nonempty factors follows immediately from
the number of blocks.

\end{proof}

  \subsection{Kernel, splitting, and base change}

With Lemma \ref{lem:classical-transpose} established, the argument of Section \ref{sec:unimodularity-divided-power} applies, with a few suitable modifications.

To simplify notation, instead of writing \((q_d^-)^*\) (resp. \((q_d^+)^*\)) in the orthogonal (resp. symplectic) case we are going to write \((q_d^{\varepsilon})^{*}\) where $\varepsilon=\pm 1$.\smallskip

We retain the admissibility conditions of Section 3,
expressed in the tableau conventions of Section \ref{sec:classical-standard-tableaux}.
Thus the corresponding tableau shapes have even rows
in the orthogonal case and even columns in the
symplectic case, with the dimensional restrictions
specified there. 
Equivalently, in the factor--size notation
$\lambda=(a_1,\ldots,a_s)$ of
Section~\ref{sec:classical-standard-tableaux},
admissibility requires $s\leq n$,
where $n$ is the
number of blocks occurring in Lemma \ref{lem:classical-transpose}.
We denote by $\mathcal T^\epsilon_{n,d}$ the set
of admissible standard tableaux of degree $d$.

\begin{lemma}\label{lem:classical-distinguished}
For every admissible standard tableau $T$, there exists an exterior
basis monomial $M$ such that
\[
(q_d^{\varepsilon})^{*}(M^{*}) = \pm T.
\]
\end{lemma}
\begin{proof}
Apply the construction of Lemma~\ref{lem:classical-transpose} to the successive factors of $T$.
\end{proof}

Then in the two cases

\begin{theorem}
\label{thm:classic-q-surj}\label{thm:split1}
\begin{align}
\label{eq:kernel1}
\operatorname{im}(q_d^\epsilon)^* =\bigoplus_{S\in  \mathcal T^\epsilon_{n,d}}\mathbb Z S;\quad \ker q_d^\epsilon&=\bigoplus_{S\in\mathcal T^\epsilon_{ d}\setminus \mathcal T^\epsilon_{n,d}}\mathbb Z S^\vee,\\
\label{eq:Div-basis-admissible1}
\Div^\epsilon_{\mathbb Z,d}:=\operatorname{im}q_d^\epsilon
 &=\bigoplus_{T\in\mathcal T^\epsilon_{n,d}}\mathbb Z\,q_d^\epsilon(T^\vee).
\end{align}
The image is a   direct summand of $E_{2d}$.

\end{theorem}
 \begin{proof}The proof repeats the argument of Theorem \ref{thm:unimodularity-saturation}, using Lemmas \ref{lem:classical-transpose} and \ref{lem:classical-distinguished} in place of Lemmas \ref{lem:allcoords} and \ref{lem:distinguished} together with the corresponding  straightening Theorem.
 \end{proof}

\begin{corollary}
 
\label{cor:classical-weyl}\strut\begin{enumerate}\item $\Div^\epsilon_{\mathbb Z } =\operatorname{im}q ^\epsilon=(E ^\epsilon)^{G_\epsilon}.$
\item The kernel $I_n$ of $q ^\epsilon$ is freely generated as a $\Z$--module by the elements of
the  basis dual to that of non-admissible standard tableaux.
\item 
The space $(E_{2d}^\epsilon)^{G_\epsilon}$ of degree $2d$  invariants  has a filtration by
$\GL(U)$-submodules with associated graded  modules

$$\bigoplus_{\lambda\,\mathrm{admissible},\ |\lambda|=d}\Delta_{\tau^\epsilon(\lambda)}(U^*)$$
These filtrations, including their associated graded modules, commute with base change.
\item Let $K$ be a  field, and require $\operatorname{char}K\ne2$
in the orthogonal case. Then the base-changed maps induce
isomorphisms of graded $\GL(U_K)$-algebras
\[
 (D^-\otimes K)/(\ker q^-\otimes K)\simeq
 \bigwedge(V^-_K\otimes U_K)^{\OO(V^-_K)},
 \qquad
 (D^+\otimes K)/(\ker q^+\otimes K)\simeq
 \bigwedge(V^+_K\otimes U_K)^{\Sp(V^+_K)},
\]
\end{enumerate}
where each formula is taken in its corresponding case.\smallskip

\end{corollary}
\begin{proof} 

(1)\quad The contraction images are contained in the invariant subspaces by
construction. 

In characteristic zero, the classical FFT for exterior invariants, recalled in Part 1, identifies the image       
with the invariants, so since $\Div^\epsilon_{\mathbb Z }$ is a direct summand of  $E^\epsilon$ the claim follows.\smallskip

(2)\quad This is the kernel description in Theorem~\ref{thm:classic-q-surj}. (3)\quad This follows by dualizing formula \eqref{clams1}, as in formula \eqref{lacl}.

\smallskip

(4)\quad By Theorem~\ref{thm:classic-q-surj} the spaces $ (D^\epsilon\otimes K)/(\ker q^\epsilon\otimes K)$  have rank
$\#\mathcal T_{n,d}^\epsilon$ and remain direct summands after base
change. 

 Theorem~\ref{thm:classical-field-independence}
shows that the full invariant dimension is unchanged in every
permitted characteristic. An injective map between these spaces with
the same dimension is surjective. This proves the algebra statements.

\end{proof}

\section{The case of $\SO$}

 Consider $A:=V^a =V\otimes U$ with the standard quadratic form.   Choosing a basis for $V$, we  display the $a$ vectors in $V^a$  as the columns $x_i=(x_{1,i}, x_{2,i},\dots, x_{n,i})$ of an $n\times a$ matrix $X$.\medskip

 Let $M=(\bigwedge A)^{SO(V)}$; the group  $\Z/(2)={\rm O}(V)/\SO(V)$ acts on $M$ via the determinant, and so we can decompose 
  $M=M^0\oplus M^1$ with  $$ M^0=(\bigwedge A)^{ O(V)} \quad\text{ and}\quad M^1:=\{m\in M\mid g\cdot m=\det(g)m,\ \forall g\in O(V)\}\,. $$
 
 Given an admissible partition, that is, a partition $\lambda$ with even rows contained in the rectangle   $n\times(a-1)$ (padding by zero rows to height $n$), let $\lambda^1$ contained in the rectangle $n\times a$ be obtained from $\lambda$ by adding  one   column  of length $n$. The partition $\lambda^1$ has odd rows.
 \begin{theorem}\label{fonn1}\strut 
 \begin{enumerate}\item $M$ is generated by the $O(V)$ invariants $(\bigwedge A)^{O(V)}$ and by the skew determinants $ {\det}_c$ of the maximal minors of $X$.
\item The space   $M^1$ has a basis indexed by standard tableaux of the admissible shapes $\lambda^1$.

\item The space   $M^1$ has a canonical  filtration by $GL(U)$-modules indexed by admissible partitions $\lambda$ contained in the rectangle $n\times(a-1)$ so that the associated graded space is the direct sum of the modules $\Delta_{\lambda^1}(U) $. 
\item  The product of two skew determinants $ {\det}_c$  of two   maximal minors of $X$   is the ${\rm O}(V)$  invariant given by formula  \eqref{perm} where one is transposed.
 \end{enumerate}

\end{theorem}

This  is similar to the previous case. One first remarks that the  skew determinants $ {\det}_c$ are $\SO$ invariants.

One decomposes  the  $\SO$ invariants  as the direct sum of the $\rm O$ invariants and the space  $M^1$  which transforms under $\det(g),\ g\in {\rm O}(V)$.\smallskip

Add a vector $y_{-1}$  and multiply $M^1$ by $Y=\prod_{i=1}^n y_{i,-1}$  getting $\rm  O$ invariants.\smallskip

Arguing as in Theorem \ref{fonn} one sees that the $\rm  O$ invariants so obtained are associated to standard tableaux (with even rows) and the first column of length   $n$ is filled with $-1$.   Thus also the  second column has length $n$.\smallskip

Then this invariant is of the form  $Y [i_1,\cdots,i_n] P$  with $[i_1,\cdots,i_n]$ a skew determinant (described by the second column)  and $P$ an $O(V)$-invariant. This proves (1). The remaining statements follow the same lines as in Theorem \ref{fonn} and we leave them to the reader.

\section{A general problem} In this short section we want to sketch what happens  when we mix Grassmann   with symmetric  variables in characteristic 0. In positive characteristic we do not know much.\smallskip

This, in the case of $GL(V)$ invariants,  means to study the invariant theory of
\begin{equation}\label{lacom}
\bigwedge (V\otimes U_1)\otimes \bigwedge (V^*\otimes W_1)\otimes S(V\otimes U_2)\otimes  S(V^*\otimes W_2).
\end{equation} 
That is, vectors $v_1,\cdots,v_h; v_{h+1},\cdots, v_m$  with   the first $v_i,\ i\leq h$ fermionic and the remaining bosonic and similar conditions for the covectors $\phi_1,\cdots,\phi_k; \phi_{k+1},\cdots, \phi_p$.\smallskip

By Aronhold's method we still have in characteristic 0 the generating invariants  $t_{i,j}=\langle \phi_i\mid v_j\rangle$  except that now the $t_{i,j}$ can be both even and odd and  the second fundamental theorem does not have such a clear statement since one has to use the Cauchy formulas for all four factors of formula \eqref{lacom} (see for example \cite{ser}).

This clearly gives a complicated description via general  Clebsch--Gordan or Littlewood--Richardson formulas, although Aronhold's method still yields a description of multilinear relations.\medskip

In fact there is a general, but not well understood, theory of mixed standard tableaux, due to Grosshans, Rota, Stein in this case \cite{GrosshansRotaStein}, which may give some insight but it is beyond the scope of this paper. 

\smallskip


\begin{thebibliography}{99}
\bibitem{Ab}S. Abeasis, \emph{The  $GL(V)$--invariant ideals in  $S(S^2V) $}

Rend. Mat. (6) 13 (1980), no. 2, 235--262.

\bibitem{Abd} S. Abeasis and A. Del Fra, \emph{Young diagrams and ideals of Pfaffians},
Adv. in Math. 35 (1980), no. 2, 158--178.

\bibitem{ABW}
K. Akin, D. A. Buchsbaum and J. Weyman,
{\it Schur functors and Schur complexes.}  Adv. in Math. 44 (1982), no. 3, 207--278.




\bibitem{brundan}
J. Brundan, \emph{Dense orbits and double cosets, in {\it Algebraic groups and their representations} (Cambridge, 1997)}, 259--274, NATO Adv. Sci. Inst. Ser. C: Math. Phys. Sci., 517, Kluwer Acad. Publ., Dordrecht, 1998.

\bibitem{BrunsVetter}
W.~Bruns and U.~Vetter, \emph{Determinantal Rings},
Lecture Notes in Mathematics \textbf{1327}, Springer, 1988.
\bibitem{Clausen}
M. Clausen, Dominance orders,  \emph{Capelli operators, and straightening of bideterminants,} European J. Combin. {\bf 5} (1984), no.~3, 207--222.



\bibitem{Co} B. Conrad, O. Gabber and G. Prasad, {\it Pseudo-reductive groups}, second edition, 
New Mathematical Monographs, 26, Cambridge Univ. Press, Cambridge, 2015; Appendix C (Section C.1 on classical groups).

\bibitem{Co1}B. Conrad, Reductive group schemes, in {\it Autour des sch\'emas en groupes. Vol. I}, 93--444, Panor. Synth\`eses, 42/43, Soc. Math. France, Paris, 2014.

\bibitem{DeConciniProcesi1976} C. De~Concini and C. Procesi, \emph{A characteristic free approach to invariant theory,} Advances in Math. {\bf 21} (1976), no.~3, 330--354.


\bibitem{DePE1980} C. De Concini, D. Eisenbud and C. Procesi, \emph{Young diagrams and determinantal varieties},
Invent. Math. 56 (1980), no. 2, 129--165.

\bibitem{DEP} C. De Concini, D. Eisenbud, C. Procesi, \emph{Hodge algebras}, Ast\'erisque \textbf{91} (1982).


\bibitem {depr0}  C. De Concini,   C. Procesi, \textit{The invariant theory of matrices}, A.M.S. University Lecture Series v. \textbf{69}, 151 pp. (2017).

\bibitem {dem} M. Demazure, \emph{Donn\'ees radicielles}, Exp. 21  in M. Artin et al., {\it Sch\'emas en groupes. Fasc. 6: Expos\'es 19 \`a{} 22}, Institut des Hautes \'Etudes Scientifiques, Paris, 1965.


\bibitem{Donkin}
S.~Donkin,
\emph{Rational Representations of Algebraic Groups: Tensor Products and Filtrations},
Lecture Notes in Math. \textbf{1140}, Springer-Verlag, Berlin, 1985.
\bibitem{Fulton}
W.~Fulton,
\emph{Young Tableaux: With Applications to Representation Theory and Geometry},
London Math. Soc. Student Texts \textbf{35}, Cambridge University Press,
Cambridge, 1997.

\bibitem{GrosshansRotaStein}
Grosshans, F. D.; Rota, G.-C.; Stein, J. A.,
\emph{Invariant Theory and Superalgebras},
CBMS Regional Conference Series in Mathematics, vol.~69,
American Mathematical Society, Providence, RI, 1987.

\bibitem{Jantzen}
J.~C.~Jantzen,
\emph{Representations of Algebraic Groups}, 2nd ed.,
Math. Surveys Monogr. \textbf{107}, American Mathematical Society,
Providence, RI, 2003.

\bibitem{Mathieu}
O.~Mathieu,
\emph{Filtrations of $G$-modules},
Ann. Sci. \`Ecole Norm. Sup. (4) \textbf{23} (1990), no.~4, 625--644.
\bibitem{O'Meara}O.~T. O'Meara, {\it Introduction to quadratic forms}, reprint of the 1973 edition, 
Classics in Mathematics, Springer, Berlin, 2000.
\bibitem{ser}
A.~N. Sergeev, An analog of the classical invariant theory for Lie superalgebras. I, II, Michigan Math. J. {\bf 49} (2001), no.~1, 113--146, 147--168.
\bibitem{Wang}
J.~P. Wang,
\emph{Sheaf cohomology on $G/B$ and tensor products of Weyl modules},
J. Algebra \textbf{77} (1982), no.~1, 162--185.

\bibitem{Weyl}
H.~Weyl,
\emph{The Classical Groups: Their Invariants and Representations},
2nd ed. with supplement, Princeton University Press, Princeton, NJ, 1946.

\end{thebibliography}
\end{document}